\documentclass[10pt]{article}
\usepackage{amsmath}
\usepackage{graphicx,psfrag,epsf}
\usepackage{epstopdf}
\usepackage{enumerate}
\usepackage{natbib}
\usepackage{url} 
\usepackage{longtable}
\usepackage{tabu}
\usepackage{amsthm} 
\usepackage{amsfonts}
\usepackage{bm}
\usepackage{graphicx}
\usepackage{multirow}
\usepackage[english]{babel}
\usepackage{authblk}
\newtheorem{theorem}{Theorem}
\newtheorem{lemma}{Lemma}
\newtheorem{assumption}{Assumption}

\theoremstyle{remark}
\newtheorem{remark}{Remark}

\usepackage{lscape}
\usepackage[inline]{enumitem}
\usepackage{subcaption}
\usepackage[makeroom]{cancel}

\newcommand{\sups}[2]{{#1}^{({#2})}}

\newcommand{\la}{\langle}
\newcommand{\ra}{\rangle}
\newcommand{\mb}{\boldsymbol}
\newcommand{\N}{{\mathcal N}}

\newcommand{\diag}{\textrm{diag}}
\newcommand{\pr}{\textrm{Pr}}
\newcommand{\cov}{\textrm{cov}}
\newcommand{\E}{\mathbb{E}}
\newcommand{\diff}{\textrm{d}}

\newcommand{\specialcell}[2][c]{%
  \begin{tabular}[#1]{@{}l@{}}#2\end{tabular}}
\allowdisplaybreaks

\usepackage{geometry}
\newcommand{\blind}{1}

\begin{document}

\def\spacingset#1{\renewcommand{\baselinestretch}%
{#1}\small\normalsize} \spacingset{1}


\if1\blind
{
  \title{\bf Fast moment estimation for generalized latent Dirichlet models}
  \author[1]{Shiwen Zhao\thanks{sz63@duke.edu}}
  \author[2]{Barbara E. Engelhardt\thanks{bee@princeton.edu}}
  \author[1]{Sayan Mukherjee\thanks{sayan@stat.duke.edu}}
  \author[1]{David B. Dunson\thanks{dunson@duke.edu}}
  \affil[1]{Department of Statistical Science, Duke University}
  \affil[2]{Department of Computer Science, Princeton University}
  \maketitle
} \fi

\if0\blind
{
  \bigskip
  \bigskip
  \bigskip
  \begin{center}
    {\LARGE\bf Fast moment estimation for generalized latent Dirichlet models}
\end{center}
  \medskip
} \fi

\bigskip
\begin{abstract}
  We develop a generalized method of moments (GMM) approach for fast
  parameter estimation in a new class of Dirichlet latent variable
  models with mixed data types.  Parameter estimation via GMM has been
  demonstrated to have computational and statistical advantages over
  alternative methods, such as expectation maximization, variational
  inference, and Markov chain Monte Carlo. The key computational
  advantage of our method (MELD) is that parameter estimation does not
  require instantiation of the latent variables. Moreover, a
  representational advantage of the GMM approach is that the behavior
  of the model is agnostic to distributional assumptions of the
  observations. We derive population moment conditions after
  marginalizing out the sample-specific Dirichlet latent
  variables. The moment conditions only depend on component mean
  parameters. We illustrate the utility of our approach on simulated
  data, comparing results from MELD to alternative methods, and we show the
  promise of our approach through the application of MELD to several data sets.
\end{abstract}

\noindent%
 \vfill

\spacingset{1.45} 
\section{Introduction}
\label{sec:intro}

Many modern statistical applications require the analysis of
large-scale, heterogeneous data types including continuous,
categorical, and count variables. For example, in social science,
survey data often consist of collections of different data types
(e.g., height, gender, and age); in population genetics, researchers
are interested in analyzing genotype (integer-valued) and
heterogeneous traits (e.g., blood pressure, BMI, alcoholic drinks per
day). Often data take the form of an $n \times p$ matrix
$\bm Y = (\bm y_1,\ldots,\bm y_n)^T$, with
$\bm y_i = (y_{i1},\ldots,y_{ip})^T$ a $p$ dimensional vector of
measurements of varying types for subject $i$, for $i = 1,\ldots,n$
subjects.

This paper focuses on 1) developing a new class of generalized latent
variable models for mixed data types, and 2) performing fast and
robust parameter estimation using the generalized method of moments
(GMM). Many models have been proposed to address the first goal
(reviewed below), but often these models lack robustness both
statistically and computationally.  In addition, parameter estimation
in these models is often both inefficient and unstable, which
necessitates our second goal.  Hence, there is a clear need for new
classes of models and corresponding robust and efficient approaches
for routine inference under these models in general 
applications.

For modeling mixed scale data, there are two general strategies that
have been most commonly employed in the literature.  The first is to
assume an underlying Gaussian model, and then to characterize the
dependence through a structured model of the Gaussian covariance
\citep{muthen_general_1984}.  Such models are routinely used in the
social science literature, focusing almost entirely on the case in
which data are categorical or continuous, with the categorical
variables arising via thresholding of Gaussian variables.  These
models are related to Gaussian copula models, which have been
developed for broad classes of mixed scale data also including counts
\citep{murray_bayesian_2013}.  The second approach is to define an
exponential family distribution for each of the $p$ variables,
inducing dependence between variables through generalized linear
models containing shared latent variables \citep{sammel_latent_1997,
  moustaki_generalized_2000, dunson_dynamic_2003}.

There are a number of disadvantages to the above approaches.  
The underlying Gaussian framework is restrictive
in forcing continuous variables to be Gaussian and non-continuous
variables to be categorical.  The exponential family latent variable
models are potentially more flexible, but in practice suffer in
several respects.  Most notably, there is a fundamental lack of
robustness due to the model structure, which arises because of the
dual role of the latent variables in controlling sample dependence and
the shape of the marginal distributions.  For example, if there are
two count variables, $y_{i1}$ and $y_{i2}$, that are highly
correlated, then both counts must load strongly on a similar set of
latent variables, which will induce over-dispersion in the marginal
distributions with the variance much larger than the mean. There is no
possibility of having highly correlated counts that do not have high
over-dispersion.  Copula models solve this problem by imposing a
restrictive Gaussian copula covariance structure.  Another fundamental
issue is computation, which tends to rely on expectation-maximization
(EM) or Markov chain Monte Carlo (MCMC) algorithms that alternate
between updating latent variables and population parameters,
intrinsically leading to slow convergence, inefficiency, and
instability.

One promising approach to address the above issues is to rely on a
generalized method of moments (GMM) estimator, which relates to the
second goal in this paper. The GMM estimator is robust to
misspecification of the higher-order moments, while also leading to
optimal efficiency among all estimators using only information on the
initial moments \citep{hansen_large_1982}. Applications of GMM to
structural models with latent variables have a long history, with
early examples including \citet{bentler_contributions_1983} and
\citet{anderson_structural_1988}, among others. More recently,
\citet{gallant_generalized_2013} applied GMM to a specific class of
latent variable model by defining moment conditions based on the
complete data, including the latent variables. In contrast,
\citet{bollen_model-implied_2014} relied on a model-implied
instrumental variable GMM estimator.  Generally, current GMM methods
focus on latent variable models that satisfy restrictive assumptions
or require the instantiation of latent variables in a computationally
intensive estimation algorithm.

The focus of our work is on a broad new class of latent variable
models for mixed scale data. This eliminates some of the problems of
current exponential family latent variable models, while also enabling
derivation of an efficient GMM implementation that marginalizes out
latent variables.  With the first goal in mind, we focus on generalized
mixed membership models, which incorporate {\em Dirichlet latent
  variables}.  Mixed membership models have a rich literature with
applications ranging from inference of ancestral populations in
genomic data \citep{pritchard_inference_2000,
  pritchard_association_2000} to topic modeling of documents
\citep{blei_latent_2003}. Using Dirichlet latent variables to define
cluster membership allows samples to {\em partially} belong to each of
$k$ latent components. Our model class includes latent Dirichlet
allocation \citep{blei_latent_2003} and simplex factor models
\citep{bhattacharya_simplex_2012} as special cases; however, we go beyond
 current approaches by allowing mixed scale data.

For the second goal, we develop an efficient moment tensor approach
for parameter estimation in mixed membership models with the Dirichlet
latent variables marginalized out. This objective is related to recent
moment tensor methods developed for latent variable models including
mixtures of Gaussians, hidden Markov models, mixed membership models,
and stochastic block models \citep{Arora:1439928,
  anandkumar_method_2012, anandkumar_spectral_2012, hsu_learning_2013,
  anandkumar_tensor_2014a, anandkumar_tensor_2014b}. This recent work
finds moment tensors of these latent variable models that have a
symmetric PARAFAC (parallel factors) tensor decomposition
\citep{kiers_towards_2000}. Parameter estimation using a moment tensor
approach is performed as follows: 1) third order empirical moment
tensors are transformed to an orthogonal decomposable form; 2)
orthogonal decompositions are performed using tensor power methods or
singular value decompositions; 3) parameter estimates are recovered
from the decomposition of the moment tensors. The moment tensor
approach offers substantial computational advantages over other
approximation methods, as illustrated in various
applications~\citep{tung_spectral_2014, anandkumar_tensor_2014a,
  colombo_fastmotif:_2015}.

In this work, we develop a moment tensor approach for generalized
mixed membership models. Here, parameter estimates are performed in
the GMM framework by minimizing quadratic forms of the moment
functions, and asymptotic efficiencies are analyzed.  Our approach,
Moment Estimation for Latent Dirichlet models (MELD), differs from
previous approaches in that the moment functions are defined for
second and third order moments for variables with different generative
distributions. This goes beyond previous moment tensor approaches that
construct moment tensors for homogeneous data distributions.

In Section 2, we characterize the class of generalized Dirichlet latent variable models in MELD. In
Section 3, we briefly review GMM and introduce the estimation
procedure used in MELD. Asymptotic properties of the GMM estimator are
also discussed. In Section 4, we evaluate the performance in a
simulation study. We apply our method to two examples in Section 5 and
conclude with a discussion in Section 6.

\section{Generalized Dirichlet latent variable models}
\label{model}

In this section, we specify a generalized Dirichlet latent variable
model. Let $\bm y_i = (y_{i1},\ldots,y_{ip})^T$ denote a vector of $p$
measurements having heterogeneous data types (e.g., continuous,
categorical, and counts) over subjects $i=1,\ldots,n$.  We are
interested in flexible models for joint distributions of the data. A
simple way to achieve this is to assume the elements of $\bm y_i$ are
conditionally independent given a component index
$s_i \in \{1,\ldots,k\}$, with the densities for each variable in each
component having a parametric form.  This is inflexible due to the
assumption that each subject $i$ belongs to exactly one component. We
instead allow subjects to partially belong to each of the $k$
components. We assume that the elements of $\bm y_i$ are conditionally
independent given a latent Dirichlet vector
$\bm x_i = (x_{i1},\ldots,x_{ik})^T \in \Delta^{k-1}$, with
$\Delta^{k-1}$ denoting the $(k-1)$ probability simplex. In
particular, we let
\begin{align}
  y_{ij} \sim \sum_{h=1}^k x_{ih} g_j(\bm \phi_{jh}), \label{eq:model}
\end{align}
where $g_j( \bm \phi_{jh})$ is the density of the $j$th variable
specific to component $h$.  The elements of $\bm x_i$ are
interpretable as probability weights for subject $i$ on each of $k$
components; in this setting, {\em pure} subjects have weight vectors
with all zeros except for a single one.

The corresponding density $g_j(\bm \phi_{jh})$ is absolutely
continuous with respect to a dominating measure
$(\Omega, \mathcal{H},\mu)$; in general, it can take any parametric
form, and may be constrained to belong to the exponential family. Our
model in Equation \eqref{eq:model} generalizes the traditional mixed
membership model to allow mixed data types. The parameters
$\bm \Phi = \{\bm \phi_{jh}\}_{1 \le j \le p, 1 \le h \le k}$ are
mixture component parameters shared by all subjects.  A full
likelihood specification is completed by choosing a population
distribution for the latent variable vector, $\bm x_i \sim P$, with
$P$ a distribution on the simplex $\Delta^{k-1}$.  In this work, we put
the Dirichlet distribution on this latent variable vector
$\bm x_i \sim \mbox{Dir}( \bm \alpha )$ with parameter $\bm \alpha =
(\alpha_1,\ldots,\alpha_k)^T$.

Let $\bm m_i = (m_{i1},\ldots m_{ip})^T$ denote a membership vector
for subject $i$, where $m_{ij} \in \{1,\ldots,k\}$ indicates the
component that feature $j$ in subject $i$ is generated from. We
specify the following generative model
\begin{align}
  y_{ij} \mid m_{ij} = h &\sim g_j(\bm \phi_{jh}),\quad 
  m_{ij} \mid \bm x_i \sim \mbox{Multi}(\bm x_i),\quad 
  \bm x_i \sim \mbox{Dir}(\bm \alpha). \label{eq:genModel}
\end{align}

When $\bm y_i$ is a multivariate categorical
vector and $g_j$ represents a multinomial distribution,
\citet{bhattacharya_simplex_2012} proposed a simplex factor model.
Their model reduces to latent Dirichlet allocation
\citep{blei_latent_2003} in the special case in which $y_{ij}$ is the
number of occurrences of word $j$ in document $i$, and
$g_j(\bm \phi_{jh})$ is a Poisson distribution, with scalar
$\bm \phi_{jh}$ indicating the rate of occurrence of word $j$ in topic
$h$.  We also assume $\bm x_i \sim \mbox{Dir}( \bm \alpha )$, but
allow a general form for the component specific densities.  This leads
to the following likelihood after marginalizing out the latent
variables:
\begin{align}
  p(\bm Y \mid \bm \alpha, \bm \Phi) =
  \prod_{i=1}^n \bigg[ \int \prod_{j=1}^p \bigg(\sum_{h=1}^k x_{ih}
  g_j(y_{ij} \mid \bm \phi_{jh}) \bigg) \diff P(\bm x_i) \bigg]. \notag
\end{align}

\paragraph{Categorical data:}
Let $y_{ij} \in \{1,\ldots,d_j\}$, for $j=1,\ldots,p$, so that we have
multivariate categorical data that can be organized as a $p$-way
contingency table.  Let $\mb c = (c_1,\ldots,c_p)^T$ with
$c_j \in \{1,\ldots,d_j\}$ and
$\pi_{\mb c} = \pr(\mb y_i = \mb c) = \pr(y_{i1} = c_1,\ldots,y_{ip} =
c_p)$.
Then $\mb \pi = \{\pi_{\mb c}\}$ is a probability tensor
\citep{dunson_nonparametric_2009} satisfying
  $\pi_{\mb c} \ge 0$ and 
  $\sum_{\mb c} \pi_{\mb c} = 1$. 
In this case,
$g_j(\bm \phi_{jh}) \stackrel{d}{=} \mbox{Multi}(\bm
\phi_{jh})$
is a multinomial distribution, with
$ \bm \phi_{jh} = (\phi_{jh1}, \ldots, \phi_{jhd_j})^T$
$\in \Delta^{d_j-1}$
a probability vector. As shown by \citet{bhattacharya_simplex_2012},
marginalizing out the mixture proportion $\bm x_i$ leads to a Tucker
decomposition of $\bm \pi$
\begin{align}
  \pi_{\mb c} = \sum_{h_1 = 1}^k \ldots \sum_{h_p = 1}^k
  g_{h_1,\ldots,h_p} \prod_{j=1}^p \phi_{j h_j
    c_j}. \notag
\end{align}

\paragraph{Mixed data types:}
When $\bm y_i$ includes variables of mixed
data types, $g_j(\bm \phi_{jh})$ belongs to the exponential
family of distributions with mean parameter
$\bm \phi_{jh}$.  We are interested in
estimating mean parameters
$\bm \Phi = (\bm \Phi_1; \ldots; \bm \Phi_p)$ with
$\bm \Phi_j = (\bm \phi_{j1},\ldots, \bm \phi_{jk})$. For categorical
variables $y_{ij}$, $\bm \phi_{jh}$ is a probability vector indicating
the distribution of the categories in component $h$; for
non-categorical variables $y_{ij}$, $\bm \phi_{jh}$ is a scalar term
indicating the mean parameter of that feature for component $h$.

\section{Generalized methods of moments estimation}

In this section, we provide a procedure for parameter estimation in
generalized Dirichlet latent variable models based on a GMM
framework. We first state the moment functions and then propose a 
two-stage procedure for parameter estimation.


\subsection{Moment functions used in MELD}


The idea behind method of moments (MM) algorithms is to derive a list
of moment equations that have expected value zero at the true
parameter values. Given observations $\{\bm y_i\}_{1\le i \le n}$, MM
algorithms specify $\ell$ moment functions and define a moment vector
$f({\bm y_i, \bm \theta}) = (f_1(\bm y_i,\bm \theta),\ldots, f_\ell
(\bm y_i, \bm \theta))^T$
such that $\mathbb{E}[f(\bm y_i, \bm \theta)] = {\bf 0}$ at the true
parameter
$\bm \theta = \bm \theta_0 \in \bm \Theta \subseteq
\mathbb{R}^p$.
The moment vector is then used to specify an estimator by
solving the following $\ell$ sample equations with $p$ unknowns:
\begin{align}
  \widehat{\bm \theta} \mbox{ such that } {\bm 0} = \Big[ {f}_n({\bm
  \theta}) \equiv \frac{1}{n} \sum_{i=1}^n f (\bm y_i, \bm \theta)
  \Big] \mbox{ at } {\bm \theta} = \widehat {\bm
  \theta}, \label{eq:MM} 
\end{align}
where $f_n(\bm \theta)$ is the sample estimate of
$\mathbb{E}[f(\bm y_i, \bm \theta)]$. When $f({\bm y_i},{\bm \theta})$
is linear in $\bm \theta$ (e.g., linear regression with instrumental
variables) with independent moment functions and $\ell = p$, then
$\bm \theta$ can be uniquely determined.  When $\ell > p$, the system
in Equation \eqref{eq:MM} might be over-determined, in which case
standard MM cannot be applied.

Generalized method of moments (GMM) minimizes the following
quadratic form rather than solving a linear system:  
\begin{align}
  \bm {\widehat \theta} \equiv \arg \min_{\bm \theta} \big[ Q_n(\bm
  \theta; \bm A_n) = {f}_n(\bm \theta)^T \bm A_n \,
  {f}_n(\bm \theta) \big],  \label{eq:GMM}
\end{align}
with $\bm A_n$ the asymptotically optimal positive semidefinite matrix: 
\begin{align}
  \bm A_n^{-1} = \bm S_n = \textrm{Var}[n^{1/2} {f}_n
  (\bm \theta)]. \notag
\end{align}
In the over-determined setting, minimizing the quadratic form is still well
defined.  Hansen's two stage procedure~\citep{hansen_large_1982} obtains
an initial estimate of $\bm {\widehat \theta}$ using a suboptimal weight matrix,
such as the identity.  This initial estimate is used to calculate $\bm A_n$, which 
is then plugged into Equation \eqref{eq:GMM} to obtain the final estimate~\citep{hall_generalized_2005} 
which is consistent and asymptotically normal  \citep{hansen_large_1982}.

Applying GMM to all of the parameters in the generalized Dirichlet
latent variable model is not feasible due to the 
dimensionality of the parameter space. Higher-order moment functions
are complex and involve large numbers of unknowns. We will make use of
a tensor decomposition and define moment functions that only depend on
lower order moments. This extends recent moment tensor approaches to
our proposed broad class of latent variable models
\citep{Arora:1439928, anandkumar_method_2012,
  anandkumar_spectral_2012, hsu_learning_2013,
  anandkumar_tensor_2014a, anandkumar_tensor_2014b}. One key
innovation in our procedure is that we allow for lower order
interaction moments across variables, where each variable can have a
different distribution. The moment functions will have rank one
decompositions of the component mean parameters over the heterogeneous
variables. These heterogeneous low-order moment functions lead to fast
and robust parameter estimation in MELD, as we show in the simulation
and application sections.

We now define the moment functions. Recall that $y_{ij}$ is the $j$th
variable for subject $i$, which may include various data types,
including continuous, categorical, or count data. Each categorical
variable $y_{ij}$ is encoded as a vector $\bm b_{ij}$ of length $d_j$ containing
all zeros except for a one in position $y_{ij}$.  If $y_{ij}$
is a non-categorical variable, then $\bm b_{ij} = y_{ij}$ is a scalar
value. In the following, we assume that there are
$k$ latent components and that the latent variable $\bm x_i$
follows a Dirichlet distribution with parameters
${\bm \alpha} =(\alpha_1,...,\alpha_k)^T$ and
$\alpha_0 = \sum_{h=1}^k \alpha_h$. Given an observation $\bm b_{ij}$
and hyperparameter $\bm \alpha$, 
we define moment functions based on
second and third moments:
\begin{align}
     \sups{\bm F}{2}_{jt}&(\bm y_i, \bm \Phi) = \bm b_{ij} \circ \bm b_{it}
                  -  \frac{\alpha_0}{\alpha_0 + 1} \bm \mu_j \circ \bm
                           \mu_t - \bm \Phi_j \sups{\bm \Lambda}{2}
                           \bm \Phi_t^T,  \quad 1\le j,t \le p, \, j
                           \ne t \label{eq:f2jt} \\
      \sups{\bm F}{3}_{jst}&(\bm y_i, \bm \Phi) = \bm b_{ij} \circ \bm b_{is} \circ
  \bm b_{it} \notag \\
  & - \frac{\alpha_0}{\alpha_0+2}\bigg(\bm b_{ij}
    \circ \bm b_{is} \circ \bm \mu_t  + \bm \mu_j
    \circ \bm b_{is} \circ \bm b_{it} + \bm b_{ij}
    \circ \bm \mu_s \circ \bm b_{it} \bigg) \notag \\
  & + \frac{2\alpha_0^2}{(\alpha_0+1)(\alpha_0+2)} \bm \mu_j
    \circ \bm \mu_s \circ \bm \mu_t
    - \sups{\bm \Lambda}{3} \times_1 \bm \Phi_j \times_2 \bm \Phi_s \times_3 \bm
                          \Phi_t,  \label{eq:f3jst} \\
                          &  \quad 1\le j,s,t\le p, \, j \ne s \ne t. \nonumber
\end{align}
We use $\circ$ to denote an outer product, and $\times_s$ indicates
multiplication of a tensor with a matrix for mode $s$. The means
${\bm \mu_j}$ take the values
$\bm \mu_j = \E(\bm b_{ij} \mid {\bm \Phi}) =
\frac{1}{\alpha_0}\bm \Phi_j \bm \alpha$
for $j = 1,...,p$.  The second moment function
$\sups{\bm F}{2}_{jt}(\bm y_i, \bm \Phi)$ is a $d_j \times d_t$
matrix, and the third moment function
$\sups{\bm F}{3}_{jst}(\bm y_i, \bm \Phi)$ is a
$d_j \times d_s \times d_t$ tensor. We set $d_j = 1$ when variable $j$
is non-categorical. The matrix $\sups{\bm \Lambda}{2}$ is a
$k \times k$ matrix with ${\bm \alpha}/[\alpha_0 (\alpha_0+1)]$ along
the diagonals and zero for all other entries. $\sups{\bm \Lambda}{3}$ is a three-way tensor with the
diagonal entry
$\sups{\lambda}{3}_{h} = 2 \alpha_h /[\alpha_0 (\alpha_0 +
1)(\alpha_0 + 2)]$ 
for $h = 1,\ldots,k$ and all the other entries are zero.

The following theorem states that, at the true parameter value, the
expectations of the moment functions are zero. The proof can be found in
Appendix \ref{app:theorem1}.

\begin{theorem}[Moment conditions in MELD]
  \label{momentTheorem}
  The expectations of the second moment matrix
  $\sups{\bm F}{2}_{jt}(\bm y_i, \bm \Phi)$ and third moment tensor
  $\sups{\bm F}{3}_{jst}(\bm y_i, \bm \Phi)$ defined in Equations
  \eqref{eq:f2jt} and \eqref{eq:f3jst} are zero at true model
  parameter values $\bm \Phi_0$
\begin{equation*} \mathbb{E}[\sups{\bm F}{2}_{jt}(\bm y_i, \bm \Phi_0)] = {\bm 0},
\quad \mathbb{E}[\sups{\bm F}{3}_{jst}(\bm y_i, \bm \Phi_0)] = {\bm 0},
\end{equation*}
with the expectations taken with respect to $\bm y_i$.
\end{theorem}

\subsection{Two stage optimal estimation}

We will state two versions of Hansen's two stage optimal GMM
estimation procedure for parameter inference in MELD, with the first
version using the second moment matrix in \eqref{eq:f2jt} and second
version using both the second moment matrix in \eqref{eq:f2jt} and the
third moment tensor in \eqref{eq:f3jst}. To fix notation, we
re-state Hansen's two stage GMM procedure as:
\begin{enumerate}
\item[(1)] estimate
  $ \widehat {\bm \theta} = \arg \min_{\bm \theta} \big[ Q_n(\bm
  \theta; \bm I) = {f}_n(\bm \theta)^T {f}_n(\bm
  \theta) \big];$
\item[(2)] given $\widehat {\bm \theta}$ calculate $\bm S_n$ and set
  $\bm A_n = \bm S_n^{-1}$;
\item[(3)] compute
  $\widehat{\bm \theta} = \arg \min_{\bm \theta} \big[ Q_n(\bm \theta;
  \bm A_n) = {f}_n(\bm \theta)^T \bm A_n {f}_n(\bm
  \theta) \big]$ as the final parameter estimate.
\end{enumerate}

We now define two moment vectors used in MELD. We define the first
version of moment vector $\sups{\bm f}{2}(\bm y_i,\bm \Phi)$ by
stacking second moment matrices in Equation \eqref{eq:f2jt} as follows:
\begin{align}
\sups{\bm f}{2}(\bm y_i,\bm \Phi) = \bigg(&\mbox{vec}[\sups{\bm
  F}{2}_{12}(\bm y_i,\bm \Phi)]^T, \ldots, \mbox{vec}[\sups{\bm
  F}{2}_{1p}(\bm y_i,\bm \Phi)]^T, \notag \\
  & \mbox{vec}[\sups{\bm F}{2}_{23}(\bm y_i,\bm \Phi)]^T, \ldots,
  \mbox{vec}[\sups{\bm F}{2}_{2p}(\bm y_i,\bm \Phi)]^T, \ldots, \mbox{vec}[\sups{\bm
  F}{2}_{p-1,p}(\bm y_i,\bm \Phi)]^T\bigg)^T. \label{eq:f2}
\end{align}
As the moment matrix is symmetric in the sense that
$\sups{\bm F}{2}_{jt}(\bm y_i,\bm \Phi) = \big[\sups{\bm
  F}{2}_{tj}(\bm y_i,\bm \Phi)\big]^T$,
we only consider the case where $j<t$. 
Assuming $d$ levels for each element of $\bm y_i$ results in a moment vector
of dimension $p(p-1)d^2/2$. The second version of moment vector
$\sups{\bm f}{3}(\bm y_i,\bm \Phi)$ is defined by stacking the
second moment matrices and third moment tensors in Equations
\eqref{eq:f2jt} and \eqref{eq:f3jst}
\begin{align}
  \sups{\bm f}{3}(\bm y_i,\bm \Phi) = \bigg(&\mbox{vec}[\sups{\bm
  F}{2}_{12}(\bm y_i,\bm \Phi)]^T, \ldots,
  \mbox{vec}[\sups{\bm F}{2}_{p-1,p}(\bm y_i,\bm \Phi)]^T, \notag \\
  & \mbox{vec}[\sups{\bm
  F}{3}_{123}(\bm y_i,\bm \Phi)]^T, \ldots, \mbox{vec}[\sups{\bm
  F}{3}_{12p}(\bm y_i,\bm \Phi)]^T, \notag \\
  & \mbox{vec}[\sups{\bm F}{3}_{134}(\bm y_i,\bm \Phi)]^T,\ldots,
    \mbox{vec}[\sups{\bm F}{3}_{13p}(\bm y_i,\bm \Phi)]^T \ldots, \mbox{vec}[\sups{\bm
  F}{3}_{p-2,p-1,p}(\bm y_i,\bm \Phi)]^T\bigg)^T.  \label{eq:f3}
\end{align}
The moment matrices $\sups{\bm F}{2}_{jt}(\bm y_i,\bm \Phi)$ follow
the same order as in $\sups{\bm f}{2}(\bm y_i,\bm \Phi)$. The
third moment tensors $\sups{\bm F}{3}_{jst}(\bm y_i,\bm \Phi)$ are
ordered such that the index on the right runs faster than the
index on the left in the subscript. The third moment tensor is
symmetric with respect to any permutations of its indices, so we include only the moment tensors with $j<s<t$. The
dimension of this moment vector is
$p(p-1) d^2/2 + \left[p^3 - 3p(p-1) - p \right]d^3/6$ with the
assumption of $d$ levels for each element of $\bm y_i$.

Based on the two moment vector versions, we define the following two
quadratic functions used for our GMM estimation:
\begin{align}
\sups{Q}{2}_n(\bm \Phi; \sups{\bm A}{2}_n) &= \sups{\bm f}{2}_n(\bm
  \Phi)^T \sups{\bm A}{2}_n \sups{\bm f}{2}_n(\bm
  \Phi), \label{eq:GMMQ2_n} \\
\sups{Q}{3}_n(\bm \Phi; \sups{\bm A}{3}_n) &= \sups{\bm f}{3}_n(\bm
  \Phi)^T \sups{\bm A}{3}_n \sups{\bm f}{3}_n(\bm
  \Phi), \label{eq:GMMQ3_n}
\end{align}
where $\sups{\bm f}{2}_n(\bm \Phi)$ and
$\sups{\bm f}{3}_n(\bm \Phi)$ are sample estimates of the
expectations of the moment vectors
\begin{align}
  \sups{\bm f}{2}_n(\bm \Phi) = \frac{1}{n}\sum_{i=1}^n \sups{\bm
  f}{2}(\bm y_i,\bm \Phi), \quad \sups{\bm f}{3}_n(\bm \Phi)
  = \frac{1}{n}\sum_{i=1}^n \sups{\bm 
  f}{3}(\bm y_i,\bm \Phi), \notag
\end{align}
and $\sups{\bm A}{2}_n$ and $\sups{\bm A}{3}_n$ are two positive
semidefinite matrices.
We emphasize that when calculating
$\sups{\bm f}{2}_n(\bm \Phi)$, the $\bm \mu_j$ and
$\bm \mu_t$ in $\sups{\bm f}{2}(\bm y_i,\bm \Phi)$ are replaced by
their sample estimates
  $\widehat{\bm \mu}_j = \frac{1}{n}\sum_{i=1}^n \bm b_{ij}$ and 
  $\widehat{\bm \mu}_t = \frac{1}{n}\sum_{i=1}^n \bm b_{it}$, 
instead of their parametric counterparts
$\frac{1}{\alpha_0} \bm \Phi_j \bm \alpha$ and
$\frac{1}{\alpha_0}\bm \Phi_t \bm \alpha$. 
In doing this the second moment matrix
$\sups{\bm F}{2}_{jt}(\bm y_i, \bm \Phi)$ becomes a function of
$\bm \Phi$ only through
$\bm \Phi_j \sups{\bm \Lambda}{2} \bm \Phi^T_t$, which is a rank one
summation of $\bm \phi_{jh}$ and $\bm \phi_{th}$. This method
facilitates development of a fast coordinate descent algorithm for
parameter estimation. Similar methods are used for calculating
$\sups{\bm f}{3}_n(\bm \Phi)$.

In the first stage, we set $\sups{\bm A}{\cdot}_n$ to an identity
weight matrix. Then the quadratic functions in Equations
\eqref{eq:GMMQ2_n} and \eqref{eq:GMMQ3_n} can be re-written as follows
\begin{align}
  \sups{Q}{2}_n(\bm \Phi, \bm I) &= \sum_{j=1}^{p-1} \sum_{t = j+1}^p
      ||\sups{\bm F}{2}_{n, jt}(\bm \Phi) ||^2_F, \notag \\
  \sups{Q}{3}_n(\bm \Phi, \bm I) &=  \sum_{j=1}^{p-1} \sum_{t = j+1}^p
      ||\sups{\bm F}{2}_{n,jt}(\bm \Phi)||^2_F +
                            \sum_{j=1}^{p-2} \sum_{s = j+1}^{p-1}
                            \sum^{p}_{t = s+1}
    ||\sups{\bm F}{3}_{n, jst}(\bm \Phi)||_F^2, \notag
\end{align}
where $\sups{\bm F}{2}_{n,jt}(\bm\Phi)$ and
$\sups{\bm F}{3}_{n, jst}(\bm \Phi)$ are defined as
\begin{align}
  \sups{\bm F}{2}_{n,jt}(\bm\Phi) = \frac{1}{n}\sum_{i=1}^n \sups{\bm
  F}{2}_{jt}(\bm y_i, \bm\Phi), \quad 
  \sups{\bm F}{3}_{n,jst}(\bm\Phi) = \frac{1}{n}\sum_{i=1}^n \sups{\bm
  F}{3}_{jst}(\bm y_i, \bm\Phi), \notag
\end{align}
with $\bm \mu_j$, $\bm \mu_s$ and $\bm \mu_t$ replaced by
$\widehat{\bm \mu}_j$, $\widehat{\bm \mu}_s$ and $\widehat{\bm \mu}_t$
respectively in $\sups{\bm F}{2}_{jt}(\bm y_i, \bm\Phi)$ and
$\sups{\bm F}{3}_{jst}(\bm y_i, \bm\Phi)$. $||\cdot ||_F^2$
indicates the Frobenius norm, the element-wise sum of squares.

We obtain a first stage estimator of $\bm \Phi$ by minimizing the
quadratic forms using Newton-Raphson. 
Note that after we substitute $\bm \mu_j$ by its sample estimate, only
the last term of $\sups{\bm F}{2}_{n,jt}(\bm \Phi)$ and
$\sups{\bm F}{3}_{n,jst}(\bm \Phi)$ involves unknown parameter
$\bm \Phi$. For simplicity we denote
\begin{align}
  \sups{\bm E}{2}_{n,jt} &= \sups{\bm F}{2}_{n,jt}(\bm \Phi) + \bm \Phi_j
                               \sups{\bm \Lambda}{2} \bm \Phi_t^T, \notag \\
  \sups{\bm E}{3}_{n,jst} & = \sups{\bm F}{3}_{n,jst}(\bm \Phi) + \sups{\bm
                                \Lambda}{3} \times_1 \bm \Phi_j \times_2 \bm
                                \Phi_s \times_3 \bm
                                \Phi_t. \notag
\end{align}
The two quantities $\sups{\bm E}{2}_{n,jt}$ and
$\sups{\bm E}{3}_{n,jst}$ can be computed directly from the
samples. We optimize $\bm \phi_{jh}$ with other mean parameters fixed.
After calculating the gradient and Hessian of
$\sups{Q}{2}_n(\bm \phi_{jh}, \bm I)$, the update rule simply becomes
\begin{align}
  \bm \phi_{jh}^s = \dfrac{\sum\limits_{t=1, t\ne
  j}^p (\sups{\overline{\bm
  E}}{2}_{n,jt} \bm \phi_{th} )^T}{(\sups{\lambda}{2}_h)
  \sum\limits_{t=1,t\ne j}^p \bm  \phi_{th}^T \bm
  \phi_{th}}, \label{eq:updatePhiM2}
\end{align}
where
$\sups{\overline{\bm E}}{2}_{n,jt} = \sups{\bm E}{2}_{n,jt} -
\sum_{h' \ne h} \sups{\lambda}{2}_{h'} \bm \phi_{jh'} \circ \bm
\phi_{th'}$
and $\sups{\lambda}{2}_h$ is the $h$th diagonal entry of
$\sups{\bm \Lambda}{2}$.

The update rule for $\bm \phi_{jh}$ with $\sups{Q}{3}_n(\bm \phi_{jh},
\bm I)$
can be calculated as
\begin{align}
  \bm \phi_{jh}^s = \frac{ \sups{\lambda}{2}_h \sum\limits_{t=1,t\ne
  j}^p  (\sups{\overline{\bm
  E}}{2}_{n,jt} \bm \phi_{th})^T + \sups{\lambda}{3}_h \sum\limits_{s=1,s\ne j}^p \bigg[
    \sum\limits_{t=1,t\ne s, t\ne j}^p (\sups{\overline{\bm E}}{3}_{n,jst}
    \times_2 \bm \phi_{sh} \times_3 \bm
  \phi_{th})^T\bigg]}{(\sups{\lambda}{2}_h)^2 \sum\limits_{t=1,t\ne
                j}^p\bm \phi_{th}^T \bm \phi_{th}  +
                     (\sups{\lambda}{3}_h)^2  \sum\limits_{s=1,s\ne j}^p
    \bigg[ \sum\limits_{ t=1,t\ne s, t \ne j}^p (\bm \phi^T_{sh} \bm \phi_{sh}) (\bm \phi^T_{th}
                \bm \phi_{th})\bigg]}, \label{eq:updatePhiM3}
\end{align}
where
$\sups{\overline{\bm E}}{3}_{n,jst} = \sups{\bm E}{3}_{n,jst} -
\sum_{h' \ne h} \sups{\lambda}{3}_{h'} \bm \phi_{jh'} \circ \bm
\phi_{sh'} \circ \bm \phi_{th'}$
and $\sups{\lambda}{3}_h$ is the $h$th diagonal entry in
$\sups{\bm \Lambda}{3}$.  The derivations can be found in Appendix
\ref{app:NRupdate}.  After updating $\bm \phi_{jh}$ using the above
equations, we retract $\bm \phi_{jh}$ to its probability simplex when
$y_{ij}$ is a categorical variable. We use the difference of the
objective function between two iterations divided by the dimension of
the moment vector to determine convergence. In particular, we stop
iterations when this value is smaller than $1\times 10^{-5}$.

After an initial consistent estimator of $\bm \Phi$ is found, we
calculate the asymptotic covariance matrix of the two versions of
moment vectors $\sups{\bm S}{\cdot}_n$ and define a new weight matrix
$\sups{\bm A}{\cdot}_n = (\sups{\bm S}{\cdot}_n)^{-1}$ for a second
stage GMM estimation. The form of $\sups{\bm S}{\cdot}_n$ can be
derived analytically, and we provide the results in the Supplement
Materials. In our implementation, the calculation of
$\sups{\bm A}{\cdot}_n$ requires the inversion of a full-rank dense
matrix $\sups{\bm S}{\cdot}_n$ with dimension scaling as $O(p^2d^2)$
for $\sups{\bm f}{2}_n(\bm \Phi)$ and $O(p^3d^3)$ for
$\sups{\bm f}{3}_n(\bm \Phi)$.  In addition, when including
the off-diagonal terms in the weight matrix, the updating rules become
intrinsically complicated.  In practice, we only extract the diagonal
elements of $\sups{\bm S}{\cdot}_n$ and let
$\sups{\bm A}{\cdot}_n = 1/\mbox{diag}[(\sups{\bm S}{\cdot}_n)]$ in
the second stage estimation. This approximation has been used in
previous GMM implementations \citep{joreskog_new_1987}. The gradient
descent update equations can be found by slight modification of
Equations \eqref{eq:updatePhiM2} and \eqref{eq:updatePhiM3} with
weights included. 

Note that the moment functions do not solve the identifiability
problems with respect to $\bm \Phi$.  When $\bm \alpha$ is a constant
vector, any permutation $\tau$ of $1,\ldots,k$ with
$\bm \Phi_j(\tau) = (\bm \phi_{j\tau(1)},\ldots, \bm \phi_{j\tau(k)})$
for all $j = 1,\ldots,p$ satisfies the moment condition. This problem
is inherited from the label switching problem in mixture models. A
similar situation occurs when there are ties in $\bm \alpha$ and the
permutation is restricted to each tie. However, in real applications,
a minimizer of the quadratic function is generally sufficient to find
a parameter estimate that is close to a unique configuration of the
true parameter. 

\subsection{Properties of parameter estimates}

We use GMM asymptotic theory to show that parameter estimates in MELD
are consistent. We assume the following regularity conditions on the
two versions of moment vector $\sups{\bm f}{\cdot}(\bm y_i, \bm \Phi)$
and the parameter space $\bm \Theta$.
\begin{assumption}[Regularity conditions (Assumption 3.2, 3.9 and
  3.10 \citep{hall_generalized_2005})]
  \label{assumption1}
  \begin{enumerate*}
  \item[1)] $\sups{\bm f}{\cdot}(\bm y_i, \bm \Phi)$ is continuous on
    $\bm \Theta$ for all $\bm y_i \in \mathcal{Y}$;
  \item[2)] $\E[\sups{\bm f}{\cdot}(\bm y_i, \bm \Phi)] < \infty$ and
    continuous for $\bm \Phi \in \bm \Theta$;
  \item[3)] $\bm \Theta$ is compact and
    $\E[\mbox{sup}_{\bm \Phi \in \bm \Theta} ||\sups{\bm f}{\cdot}(\bm y_i,
    \bm \Phi)||] < \infty$.
  \end{enumerate*}
\end{assumption}
\begin{remark}
  Conditions 1) and 2) are satisfied in MELD. Condition 3) is
  satisfied automatically for categorical variables, noting that
  $\bm \phi_{jh} \in \Delta^{d_j-1} \subset \mathbb{R}^{d_j}$ is
  compact. For non-categorical variables, we can restrict support to a
  large compact domain such as closed intervals in $\mathbb{R}$
  without sacrificing practical performance.
\end{remark}

With these conditions, we further assume that the weight matrix
$\sups{\bm A}{\cdot}_n$ converges to a positive definite matrix
$\sups{\bm A}{\cdot}$ in probability.  We define the population
analogs of the quadratic functions as
\begin{align}
  \sups{Q}{2}_0(\bm \Phi; \sups{\bm A}{2} )
  &=  \E[\sups{\bm f}{2}(\bm y_i,\bm \Phi)]^T \sups{\bm A}{2}
    \E[\sups{\bm f}{2}(\bm y_i,\bm \Phi)], \label{eq:GMMQ2_0} \\
  \sups{Q}{3}_0(\bm \Phi; \sups{\bm A}{3} )
  &=  \E[\sups{\bm f}{3}(\bm y_i,\bm \Phi)]^T \sups{\bm A}{3}
    \E[\sups{\bm f}{3}(\bm y_i,\bm \Phi)]. \label{eq:GMMQ3_0}
\end{align}

We have the following lemma showing the uniform convergence of
$\sups{Q}{\cdot}_n(\bm \Phi; \sups{\bm A}{\cdot}_n)$.

\begin{lemma}[Uniform convergence (Lemma 3.1
  \citep{hall_generalized_2005})]
  \label{uniformLemma}
  Under regularity conditions in Assumption \ref{assumption1},
  \begin{align}
    \sup_{\bm \Phi \in \bm \Theta}|\sups{Q}{2}_n(\bm \Phi; \sups{\bm A}{2}_n) -
    \sups{Q}{2}_0(\bm \Phi; \sups{\bm A}{2})| \overset{p}{\rightarrow} 0 \qquad
    \sup_{\bm \Phi \in \bm \Theta}|\sups{Q}{3}_n(\bm \Phi; \sups{\bm A}{3}_n) -
    \sups{Q}{3}_0(\bm \Phi; \sups{\bm A}{3})| \overset{p}{\rightarrow} 0. \notag
  \end{align}
\end{lemma}

\begin{theorem}[Consistency]
  \label{consistencyTheorem}
  Under the same conditions in Lemma \ref{uniformLemma}, the estimator
  $\sups{\widehat{\bm \Phi}}{2}$ that minimizes
  $\sups{Q}{2}_n(\bm \Phi; \sups{\bm A}{2}_n)$ converges to the true
  parameter $\bm \Phi_0$ in probability. A similar result holds for
  $\sups{\widehat{\bm \Phi}}{3}$ that minimizes
  $\sups{Q}{3}_n(\bm \Phi; \sups{\bm A}{3}_n)$.
\end{theorem}
The proof can be found in Appendix \ref{app:theorem2}. 

The asymptotic normality of $\sups{\widehat{\bm \Phi}}{2}$ and
$\sups{\widehat{\bm \Phi}}{3}$ can also be established by assuming
the following conditions on
$\partial \sups{\bm f}{\cdot}(\bm y_i, \bm \Phi) / \partial \bm \Phi$.
\begin{assumption}[Conditions on
  $\partial \sups{\bm f}{\cdot}(\bm y_i, \bm \Phi) / \partial \bm
  \Phi$ (Assumptions 3.5, 3.12 and 3.13. \citep{hall_generalized_2005})]
  \label{assumption2}
  \begin{enumerate*}
  \item[1)]
    $\partial \sups{\bm f}{\cdot}(\bm y_i, \bm \Phi)/ \partial
    \bm \Phi$
    exists and is continuous on $\bm \Theta$ for all
    $\bm y_i \in \mathcal{Y}$;
  \item[2)] $\bm \Phi_0$ is an interior point of $\bm \Theta$;
  \item[3)] $\E[\partial \sups{\bm f}{\cdot}(\bm y_i, \bm \Phi)/ \partial
    \bm \Phi] = \sups{\bm G}{\cdot}_0(\bm \Phi) < \infty$;
  \item[4)] $\sups{\bm G}{\cdot}_0(\bm \Phi)$ is continuous on some neighborhood
    $N_\epsilon$ of $\bm \Phi_0$;
  \item[5)] the sample estimate $\sups{\bm G}{\cdot}_n(\bm \Phi)$ uniformly converges
    to $\sups{\bm G}{\cdot}_0(\bm \Phi)$.
  \end{enumerate*}
\end{assumption}

\begin{remark}
  We derive
  $\partial \sups{\bm f}{\cdot}(\bm y_i, \bm \Phi)/ \partial \bm \Phi$
  in Appendix \ref{app:derivatives}. Conditions 1, 3, and 4 are
  satisfied in MELD. Condition 5 can be shown with continuousness of
  the derivative and the compactness of $\bm \Theta$.
\end{remark}

\begin{theorem}[Asymptotic normality]
  \label{normalityTheorem}
  With Assumptions \ref{assumption1} and \ref{assumption2}, we have
  \begin{align}
    n^{1/2} \bigg(\mbox{vec}(\sups{\widehat{\bm \Phi}}{\cdot}) -
    \mbox{vec}(\bm \Phi_0)\bigg) \overset{p}{\rightarrow}
    \N\bigg(\bm 0, \sups{\bm M}{\cdot} \sups{\bm S}{\cdot}
    (\sups{\bm M}{\cdot})^T\bigg) \notag
  \end{align}
  with
  \begin{align}
    \sups{\bm M}{\cdot} &= [(\sups{\bm G}{\cdot}_0)^T \sups{\bm A}{\cdot}
    \sups{\bm G}{\cdot}_0]^{-1}(\sups{\bm G}{\cdot}_0)^T \sups{\bm A}{\cdot}, \notag
  \end{align}
\end{theorem}
{\noindent where
$\sups{\bm G}{\cdot}_0 = \E[\partial \sups{\bm f}{\cdot}(\bm y_i, \bm
\Phi)/ \partial \bm \Phi] |_{\bm \Phi = \bm \Phi_0}$
and
$\sups{\bm S}{\cdot} = \lim_{n\rightarrow \infty} \mbox{Var}[n^{1/2}
\sups{\bm f}{\cdot}_n(\bm \Phi_0)]$.
The proof can be found in Appendix \ref{app:theorem3}. The optimal
estimator can be obtained so that the weight matrix
$\sups{\bm A}{\cdot}_n \rightarrow \sups{\bm A}{\cdot} = (\sups{\bm
  S}{\cdot})^{-1}$ \citep{hansen_large_1982}.}

\subsection{Model selection using goodness of fit tests}

We use goodness of fit tests to choose the number of latent components
$k$. With the optimal weight matrix
$\sups{\bm A}{\cdot}_n \rightarrow (\sups{\bm S}{\cdot})^{-1}$, 
Wald-, Lagrange multiplier- (score) and
likelihood ratio-type tests can be constructed. We could construct a sequence of test
statistics with increasing $k$ to quantify the goodness of fit in
MELD.  However, this approach requires the calculation of the optimal
weight matrix and large matrix inversions. Instead, we quantify
goodness of fit using a fitness index
(FI)~\citep{bentler_contributions_1983}. Practically, we show in an
extensive simulation study that the FI has low error in estimating $k$
in MELD.


The FI is defined based on the value of the objective function
evaluated at the parameter estimate~\citep{bentler_contributions_1983}
\begin{align}
  \mbox{FI} = 1 - \frac{\sups{Q}{\cdot}_n(\sups{\widehat{\bm \Phi}}{\cdot},
  \sups{\bm A}{\cdot}_n)}{(\sups{\bm e}{\cdot}_n)^T \sups{\bm
  A}{\cdot}_n  \sups{\bm e}{\cdot}_n}, \label{eq:FI}
\end{align}
where $ \sups{\bm e}{2}_n$ is the concatenation of
$\mbox{vec}(\sups{\bm E}{2}_{n,jt})$ for $j < t$ and
$ \sups{\bm e}{3}_n$ is the concatenation of both
$\mbox{vec}(\sups{\bm E}{2}_{n,jt})$ for $j < t$ and
$\mbox{vec}(\sups{\bm E}{3}_{n,jst})$ for $j < s < t$. This FI
is for any weight matrix $\sups{\bm A}{\cdot}_n$.  It can be viewed as
a normalized objective function; thus, the FI has values less than
one. Larger values of FI indicate a better fit.

\section{Simulation study}

In this section, we evaluate the accuracy and run time of MELD in
simulations with both categorical and mixed data types. We use two
stage estimations described in previous sections. In the first stage
an identity weight matrix is used. In second stage we set
$\sups{\bm A}{\cdot}_n = 1/\mbox{diag}[(\sups{\bm S}{\cdot}_n)]$. For
notation convenience we suppress the weight matrix
$\sups{\bm A}{\cdot}$ and the subscript $n$ in the objective functions
$\sups{Q}{\cdot}_n(\bm \Phi, \sups{\bm A}{\cdot})$.

\subsection{Categorical data}

For simulation with categorical data, we considered a low dimensional
setting ($p=20$) so that both second and third order moment functions
may be efficiently calculated. A simulation with a moderate dimension
($p=100$) was also studied and is summarized in the Supplementary
Materials.

We assumed each of the $20$ variables has $d = 4$ levels. We set the
number of components to $k = 3$ and generated $\bm \phi_{jh}$ from
$\mbox{Dir}(0.5,0.5,0.5,0.5)$ with $h=1,\ldots,k$, and $\bm \alpha$
was set to $(0.1,0.1,0.1)^T$. We drew $n = \{50,100,200,500,1,000\}$
samples from the generative model in Equation \eqref{eq:model} with
$g_j(\bm \phi_{jh})$ a multinomial distribution. For each value of
$n$, we generated ten independent data sets. We ran MELD for different
values of $k = \{1,\ldots,5\}$. The FI criterion consistently chose
the correct number of latent components $k$ in first stage estimation
(Table \ref{tab:Fitness}). For second stage, FI did not perform
well. The trajectories of the objective functions under different
values of $k$ are shown in Figure S1. The convergence of parameter
estimations on the ten simulated data sets under $n=1,000$ and $k=3$
is shown in Figure S2. MELD converged in about $25$ iterations with
$\sups{Q}{2}(\bm \Phi)$ and in about $10$ iterations with
$\sups{Q}{3}(\bm \Phi)$.

\begin{table}[h!]
  \caption{{\bf Goodness of fit tests using the fitness index (FI) in
    categorical simulation.} Larger values of FI indicate
    better fit, with the maximum at one. Results
    shown are based on ten simulated data sets for each value of
    $n$. Standard deviations of FI are provided in parentheses.
 } \label{tab:Fitness}
  \tiny
  \begin{center}
    \begin{tabular}{*{2}{l}*{4}{c}}
      \hline
      $n$ & $k$ & {$\sups{Q}{2}(\bm \Phi)$ 1st stage} & {$\sups{Q}{2}(\bm \Phi)$ 2nd stage} 
      & {$\sups{Q}{3}(\bm \Phi)$ 1st stage} & {$\sups{Q}{3}(\bm \Phi)$ 2nd stage} \\
      \hline
      \multirow{5}{*}{50}
      & 1 & 0.824(0.020) & -1.865(4.950) & 0.585(0.031) & -27.987(40.128) \\
      & 2 & 0.908(0.011) & -1.285(2.630)  & 0.745(0.032) & -27.240(54.415) \\
      & \bf 3 & \bf 0.930(0.004) & 0.570(0.044) & \bf 0.775(0.022) & -27.713(42.607) \\
      & 4 & 0.901(0.010) & 0.588(0.032) & 0.735(0.023) & 0.254(0.090) \\
      & 5 & 0.795(0.033) & \bf 0.686(0.021) & 0.653(0.024) & \bf 0.305(0.091) \\
      \hline
      \multirow{5}{*}{100}
      & 1 & 0.860(0.012) & -0.521(2.142) & 0.651(0.021) & -0.031(0.703) \\
      & 2 & 0.930(0.011) & 0.282(0.644) & 0.795(0.030) & -4.457(8.924)\\
      & \bf 3 &\bf 0.960(0.005) & 0.677(0.044) & \bf 0.851(0.009) & -4.868(11.889)\\
      & 4 & 0.942(0.012) & 0.691(0.042) & 0.822(0.010) & 0.225(0.019)\\
      & 5 & 0.863(0.046) & \bf 0.782(0.022) & 0.768(0.044) & \bf 0.232(0.080)\\
      \hline
      \multirow{5}{*}{200}
      & 1 & 0.869(0.012) & 0.679(0.060) & 0.682(0.021) & 0.298(0.070)\\
      & 2 & 0.940(0.007) & 0.699(0.047) & 0.838(0.014) & \bf 0.306(0.054)\\
      & \bf 3 & \bf 0.980(0.001) & 0.761(0.061) & \bf 0.919(0.004) & 0.278(0.049) \\
      & 4 & 0.967(0.006) & 0.780(0.019) & 0.891(0.008) & 0.287(0.050)\\
      & 5 & 0.911(0.017) & \bf 0.824(0.008) & 0.864(0.015) & 0.286(0.042) \\
      \hline
      \multirow{5}{*}{500}
      & 1 & 0.882(0.007) & 0.783(0.022) & 0.713(0.012) & 0.414(0.080)\\
      & 2 & 0.948(0.006) & 0.799(0.019) & 0.870(0.013) &\bf 0.427(0.080)\\
      & \bf 3 & \bf 0.992($<$0.001) & 0.884(0.024) &\bf 0.966(0.001) & 0.388(0.073)\\
      & 4 & 0.983(0.004) & 0.874(0.026) & 0.938(0.005) & 0.365(0.065)\\
      & 5 & 0.937(0.014) & \bf 0.894(0.006) & 0.921(0.007) & 0.353(0.042)\\
      \hline
      \multirow{5}{*}{1,000}
      & 1 & 0.888(0.003) & 0.828(0.006) & 0.729(0.008) & 0.571(0.017)\\
      & 2 & 0.951(0.003) & 0.855(0.009) & 0.881(0.008) &\bf 0.615(0.030) \\
      & \bf 3 & \bf 0.996($<$0.001) & \bf 0.950(0.005) &\bf 0.982(0.001) & 0.609(0.031) \\
      & 4 & 0.989(0.002) & 0.938(0.004) & 0.961(0.003) & 0.579(0.030)\\
      & 5 & 0.953(0.010) & 0.932(0.006) & 0.951(0.006) & 0.550(0.034)\\
      \hline
    \end{tabular}
  \end{center}
\end{table}

%

\begin{figure}[h!]
  \centering
  \includegraphics[width=\textwidth]{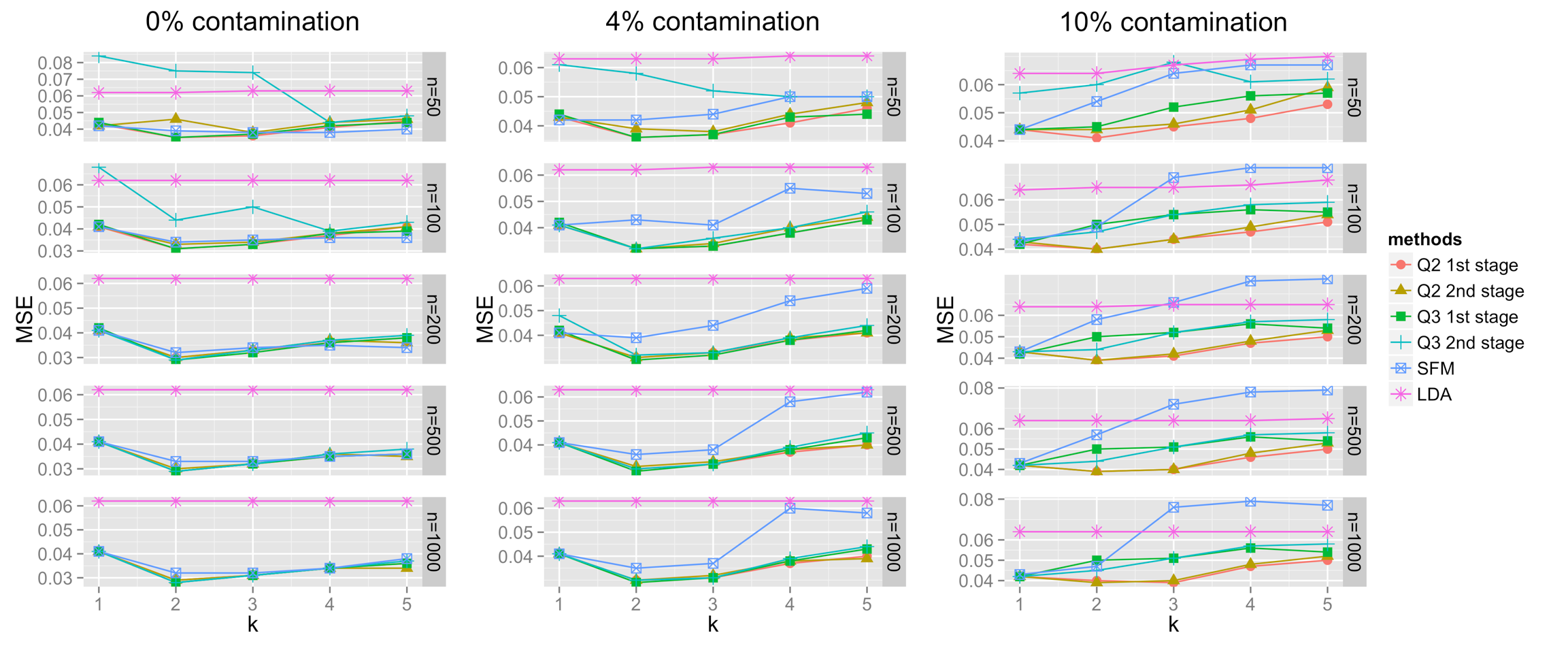}
  \caption{ {\bf Comparison of mean squared error (MSE) of estimated
      parameters in categorical simulations.} For SFM and
    LDA, posterior means of parameters are calculated using 100 posterior draws on
    each of the ten simulated data sets. The values of the MSE's and
    their standard deviations are in Supplementary Table S2,
    S3 and S4. } 
  \label{fig:CTp20MSE}
\end{figure}

We compared MELD with the simplex factor model (SFM)
\citep{bhattacharya_simplex_2012} and latent Dirichlet allocation
(LDA) \citep{blei_latent_2003}.  For SFM, we ran $10,000$ steps of
MCMC with fixed $k$ and a burn-in of $5,000$ iterations. Posterior
thinned samples were collected by keeping one posterior draw after
every $50$ steps. For the LDA model, we used the \verb$lda$ package in
R \citep{chang_lda:_2012} with collapsed Gibbs sampling. We used the
same number of MCMC iterations and burn-in iterations as with SFM. The
Dirichlet parameter for mixture proportions was set to $\alpha = 0.1$
and the Dirichlet parameter for topic distributions was set to
$\beta = 0.5$. We calculated mean squared errors (MSE's) of
different estimates as follows. We first recovered the membership
variable $m_{ij}$ and then calculated the parameter estimate of
$y_{ij}$. The MSE was calculated between the estimated parameters of
$y_{ij}$'s and their true parameters (see Supplementary Materials for
details).  MSE's with different values of $k$ are shown in Figure
\ref{fig:CTp20MSE} and the running times of different methods are
reported in Table S1. MELD $\sups{Q}{2}(\bm \Phi)$ with first stage
estimation had the most accurate parameter estimation and fastest
running speed in most cases. The second stage of MELD
$\sups{Q}{3}(\bm \Phi)$ did not perform well with small values of $n$
(i.e., $n=50$), but estimation accuracy was better when $n$ was
larger. SFM had comparable MSE's when $n$ was not large. However when
$n = 500$ and $1,000$, MELD outperformed SFM. 

We further evaluated performance in the
presence of contamination. For each simulated data set, we randomly
replaced a proportion of observations ($4\%$ and $10\%$) with
draws from a discrete uniform distribution. The MSE's under different values of $k$ are shown in Figure
\ref{fig:CTp20MSE}. With $4\%$ contamination, MELD had the most accurate
parameter estimation in almost all cases. MELD $\sups{Q}{2}(\bm \Phi)$
with first stage estimation performed the best, followed by MELD
$\sups{Q}{3}(\bm \Phi)$ with first stage estimation. The MSE's of SFM
increased.  When we increased contamination to $10\%$, MELD had the most
accurate MSE in all cases. MELD $\sups{Q}{2}(\bm \Phi)$ with first
stage estimation performed best, followed by MELD
$\sups{Q}{2}(\bm \Phi)$ with second stage estimation.  MELD
$\sups{Q}{2}(\bm \Phi)$ consistently performed better than
$\sups{Q}{3}(\bm \Phi)$, suggesting the robustness of using lower
order moments in parameter estimation. 

\subsection{Mixed data types}

We considered a simulation setting mimicking applications
in which DNA sequence variations influence a quantitative trait. An additional simulation combining categorical, 
Gaussian, and Poisson variables is summarized in the Supplementary Materials.
We generated a sequence of nucleotides $\{A,C,G,T\}$ at $50$ genetic loci 
along with a continuous or integer-valued trait, leading to $p=51$ variables.
We set $k=2$ latent components and simulated
$n=1,000$ individuals, with the first $500$ from the first component
and the last $500$ from the second component.  We chose eight loci
 $J = \{2, 4, 12, 14, 32, 34, 42, 44\}$ to be associated
with the trait.  Their multinomial parameters for each of the two
components were randomly drawn from $\mbox{Dir}(0.5,0.5,0.5,0.5)$.
The distributions for nucleotides in other loci were set to $\mbox{Multi}(0.25,0.25,0.25,0.25)$. 
Continuous traits were drawn from $\N(-3,1)$ and $\N(3,1)$, while 
count traits were drawn from $\mbox{Poisson}(5)$ and $\mbox{Poisson}(10)$, respectively for the 
two components. Ten data sets were simulated from the generative model in Equation
\eqref{eq:model}.  To assess robustness, we added contamination (e.g., through genotyping errors) by  
randomly replacing $4\%$, $10\%$ and $20\%$ of the nucleotides with values 
uniformly generated from $\{A,C,G,T\}$.

We ran MELD with first stage estimation of
$\sups{Q}{2}(\bm \Phi)$. We chose the number of components
$k = \{1,\ldots,5\}$. The fitness test indicated that FI chose the
correct value of $k$ on all ten data sets (Table S8).
For each genomic locus, we calculated its marginal frequency according
to the simulated data, and then we computed the averaged KL distance
between the estimated component distributions and the marginal
frequency as follows
\begin{align}
  \mbox{aveKL}(y_{ij}) = \frac{1}{k}\sum_{h=1}^k \sum_{c_j=1}^{d_j}
  \pr(y_{ij}=c_j|m_{ij} = h) \log\bigg(\frac{\pr(y_{ij}=c_j|m_{ij} =
  h)}{\pr(y_{ij}=c_j)}\bigg). \label{eq:aveKL}
\end{align}
A smaller averaged KL distance suggests that the component
distributions were closer to the marginal distribution, implying that
the locus frequency was not differentiated across components. The 
set $J$ corresponded exact to the eight loci with largest averaged KL
distance (Table \ref{tab:QTNsim}). 

We compared MELD with the Bayesian copula factor model \citep{murray_bayesian_2013}, 
which estimates a correlation matrix $\bm C$ between variables.
From this estimate, we computed partial correlations between the response variable (trait)
and each genetic locus \citep{hoff_extending_2007}. We ran MCMC for
the Bayesian copula factor model $10,000$ iterations with the correct
value for $k$ and a burn-in of $5,000$ iterations. Posterior samples
were collected every $50$ iterations.  We then selected genomic
locations for which their $95\%$ posterior intervals of the partial
correlation did not include zero. The resulting loci are shown in
Table \ref{tab:QTNsim}. The Bayesian copula factor model selected
nucleotides that are not in $J$ and missed locus $32$ in most cases.



\begin{table}[h!]
  \caption{{\bf Quantitative trait association simulation with $50$
      nucleotides and one response.} Nucleotides not
    in $J = \{2, 4, 12, 14, 32, 34, 42, 44\}$ are labeled by an underline
    and missing nucleotides are crossed out. Results shown are for one
    of the ten simulated data sets. The complete results can be found in
    Table S9 and S10.}
  \label{tab:QTNsim}
  \begin{center}
    \scriptsize
    \begin{tabular}{*{2}{l}|*{2}{c}}
      \hline
      Response & Contamination & $\sups{Q}{2}(\bm \Phi)$ 1st stage & Bayesian copula factor model\\
      \hline
      \multirow{4}{*}{Gaussian}& $0\%$ &  $\{2,4,12,14,32,34,42,44\}$ 
              & $\{2,4,12,14,\underline{18}, \underline{27},\cancel{32},34,42,44\}$ \\
      & $4\%$  &  $\{2,4,12,14,32,34,42,44\}$ 
              & $\{2,4,12,14,\underline{18},\underline{27},\cancel{32},34,42,44,\underline{45}\}$\\
      & $10\%$ &  $\{2,4,12,14,32,34,42,44\}$ 
              & $\{2,4,12,\cancel{14},\underline{27},\cancel{32},34,42,44,
                                             \underline{49},\underline{50}\}$ \\
      & $20\%$ &  $\{2,4,12,14,32,34,42,44\}$ 
              &$\{2,4,\underline{9},\underline{11},12,\cancel{14},\underline{20},
               \cancel{32},34,42,44\}$ \\
      \hline
      \multirow{4}{*}{Poisson}& $0\%$ &  $\{2,4,12,14,32,34,42,44\}$ 
              & $\{2,4,\underline{7},12,14,\cancel{32},34,42,44\}$ \\
      & $4\%$  &  $\{2,4,12,14,32,34,42,44\}$ 
              & $\{2,4,12,14,\cancel{32},34,42,44\}$ \\
      & $10\%$ &  $\{2,4,12,14,32,34,42,44\}$ 
              & $\{2,4,\underline{7},12,14,\underline{16},32,34,42,44\}$ \\
      & $20\%$ &  $\{2,4,12,14,32,34,42,44\}$ 
              &
                $\{2,4,\underline{7},12,\cancel{14},32,34,\underline{35},42,44\}$ \\
      \hline
    \end{tabular}
  \end{center}
\end{table}

\section{Applications of MELD}


\paragraph{Promoter sequence analysis} We applied MELD to gene
transcription promoter data available in the UCI machine learning
repository \citep{Lichman:2013}. The data include $n=106$ nucleotide
sequences $\{A,C,G,T\}$ of length $57$. The first $53$ sequences are
located in promoter regions, and the last $53$ sequences
are located in non-promoter regions.  The goal is
binary classification: using nucleotide sequence alone, we would like
to predict whether or not the sequence is in a promoter or a
non-promoter region.

We included the promoter or non-promoter status of the sequences
as an additional binary variable, giving us $p=57+1$ categorical
variables. We applied MELD $\sups{Q}{2}(\bm \Phi)$ with first stage
estimation on the full data and also on the subset of the sequences in
the promoter region and the subset of sequences in non-promoter
regions separately. We set
$k = \{1,\ldots,8\}$.  For $k=2$, MELD converged in $2.13$ seconds,
compared with SFA, which took $41.6$ seconds to perform $10,000$ MCMC
iterations. We evaluated different values of $k$ using the
goodness of fit test. FI selected two components for the full data,
two components for the promoter data, and one component for the
non-promoter data (Table S13).

We choose $k = 2$ in the following analysis. For each nucleotide position,
we calculated the averaged KL distance between the estimated component
distributions and its marginal distribution using Equation
\eqref{eq:aveKL}. An interpretation of the averaged KL distance is
that it quantifies the stratification of each nucleotide distribution
across components: a larger value of the averaged KL distance
indicates greater stratification across components, which suggests
that the nucleotide is important in defining and differentiating the
components. For the full data set, we observed approximately two peaks
of the averaged KL distance, one around the $15$th nucleotide and one
around the $42$nd nucleotide (Figure S4). The first
peak corresponds to the start of the biologically conserved region for
promoter sequences \citep{harley_analysis_1987}. For MELD applied only
to promoter sequences, this peak was reduced, suggesting that, at
approximately the $15$th nucleotide, the components all included
similarly well conserved distributions of this nucleotide. 
The estimated component membership
variables also showed the importance of nucleotides around $15$th
nucleotide (Figure S5 and S6).  For the peak around the $42$nd
nucleotide, this phenomenon was reversed: the increased averaged KL
distance remained in promoter sequences but diminished in non-promoter
sequences. One possible explanation is that this region was conserved
uniformly in the non-promoter region, but differentially conserved in
the promoter region.

\paragraph{Political-economic risk data}

In a second application, we applied MELD to political-economic risk
data \citep{quinn_bayesian_2004}, which include five proxy variables
of mixed types measured for 62 countries. The data set has been
collected and analyzed to quantify a sense of political-economic risk,
a latent quantity associated with each of the 62 countries, using a
mixed Gaussian and probit factor model \citep{quinn_bayesian_2004} and
a Bayesian copula factor model \citep{murray_bayesian_2013}. The data
are available in the \verb+MCMCpack+ package. There are three categorical variables and
two real valued variables (Table \ref{tab:PEriskExp}).

\begin{table}[h!]
  \caption{Variables in the political-economic risk data} \label{tab:PEriskExp}
  \footnotesize
  \begin{center}
    \begin{tabular}{|*{3}{l}|}
      \hline
      Variable & Type & Explanation \\
      \hline
      \verb+ind.jud+ & binary & \specialcell{An indicator variable that
                         measures the independence \\
                         of the national judiciary. This variable is equal to one if the\\
                       judiciary is judged to be independent and equal to zero otherwise.} \\
      \hline
      \verb+blk.mkt+ & real & \specialcell{Black-market premium
                              measurement. Original values are measured \\
                           as the black-market exchange rate (local
                           currency per dollar) \\
                           divided by the official exchange rate minus one.
                           \citet{quinn_bayesian_2004} \\ transformed the original data to log scale.} \\
      \hline
      \verb+lack.exp.risk+ & ordinal & \specialcell{Lack of appropriation
                                risk measurement.\\ 
                                Six levels with coding $0<1<2<3<4<5$.}
      \\
      \hline
      \verb+lack.corrup+ & ordinal & \specialcell{Lack of corruption measurement.\\
                    Six levels with coding $0<1<2<3<4<5$.}\\
      \hline
      \verb+gdp.worker+ & real &  \specialcell{Real gross domestic
                                product (GDP) per worker in 1985
                                international prices.\\ Recorded data
                                are log transformed.}
      \\
      \hline
    \end{tabular}
  \end{center}
\end{table}

We applied MELD with $k = \{1,\ldots,5\}$ using both
$\sups{Q}{2}(\bm \Phi)$ and $\sups{Q}{3}(\bm \Phi)$ with first stage
estimation to the data set. For $\sups{Q}{2}(\bm \Phi)$ with $k=3$
MELD converged in $0.10$s, and for $\sups{Q}{3}(\bm \Phi)$ with $k=3$
MELD converged in $0.45$s. The Bayesian copula factor model took
$0.91$s to complete $10,000$ MCMC iterations. The FI criterion for
$\sups{Q}{2}(\bm \Phi)$ selected $k=4$ and, for
$\sups{Q}{3}(\bm \Phi)$, selected $k=3$ (Table S14). We chose results
from $\sups{Q}{3}(\bm \Phi)$ with $k=3$ for further analysis.

The estimated component parameters for the five variables showed
distinct interpretations of the three components
(Figure~S7). We might interpret the three components
as low-risk, intermediate-risk, and high-risk political-economic
status respectively.  The first component had a high probability of
independence of the national judiciary (\verb+ind.jud+ being one) and
a low measurement of black-market premium.  The first component also
had a high probability of observing $4$th and $5$th levels in
\verb+lack.exp.risk+ and $3$rd, $4$th, and $5$th levels in
\verb+lack.corrup+. The mean of the GDP per worker was highest among
the three components. The second component had a
relatively high probability of being zero in \verb+ind.jud+ and a
large mean value of \verb+blk.mkt+. Both of lack of appropriation risk
measurement and lack of corruption measurement put higher weights on
lower category numbers ($0$, $1$ and $2$), indicating more risk and
higher levels of corruption; the GDP per worker was still high. We
might interpret this component as a society being relatively unstable
while still having a good economic forecast, meaning that GDP per
worker is high, possibly through the black market. The last component
had the least judicial independence as quantified by the probability
of \verb+ind.jud+ being zero. The black-market premium is also low, as
is the lack of risk level and lack of corruption level. The GDP per
worker is by far the lowest among the three components. We might
interpret this component as society being the most unstable with the
greatest economic risk. We found although the three components had
distinct stratification, each country was a mixture of the three
components (Figure S8).


\section{Discussion}
\label{sec:disc}

In this paper, we developed a new class of latent variable models with
Dirichlet-distributed latent variables for mixed data types. These
generalized latent Dirichlet variable models extend previous mixed
membership models such as LDA \citep{blei_latent_2003} and simplex
factor models \citep{bhattacharya_simplex_2012} to allow mixed data
types. For this class of models, we developed a fast parameter
estimation procedure using generalized method of moments. Our
procedure extends the moment tensor methods developed in recent work
\citep{anandkumar_tensor_2014b} to models with mixed data types. Our
approach does not require instantiation of latent variables. We derive
population moment conditions after marginalizing out the latent
Dirichlet variables. We demonstrated the utility of our approach using
simulations and two applications. Our results show that MELD is a
promising alternative to MCMC or EM methods for parameter estimation,
producing fast and robust parameter estimates. Since our method
depends only on certain forms of sample moments, parameter estimation
does not scale with sample size $n$ after the moment statistics have
been computed from observations. An online method to update moment
statistics when new samples arrive would allow re-estimation of the
parameters to include new observations. One limitation of our method
is that the Newton-Raphson method is of order $O(p^2)$ using second
moment functions and order $O(p^3)$ using third moment functions. One
possible approach to ensure tractability of MELD when $p$ is large is
to use stochastic gradient methods to calculate an approximate
gradient in each step.

\begin{center}
{\large\bf SUPPLEMENTARY MATERIALS}
\end{center}
We provide Supplementary Materials that include
\begin{enumerate*}
\item[a)] Derivations of asymptotic optimal weight matrices;
\item[b)] Calculation of MSE and computational complexity analysis;
\item[c)] Two additional simulations with categorical variables and mixed data types;
\item[d)] Supplementary Tables S1-S14;
\item[e)] Supplementary Figures S1-S8.
\end{enumerate*}

\appendix
\section*{Appendix}
\renewcommand{\thesubsection}{\Alph{subsection}}

\subsection{ Proof of Theorem \ref{momentTheorem}}
\label{app:theorem1}
\begin{proof}

  We start with the case where $y_{ij}$ is a categorical data with
  $d_j$ different levels.  The latent probability vector
  $\mb x_i = (x_{i1}, \cdots, x_{ik})^T \in \Delta^{k-1}$ defines the
  mixture proportion of individual $i$.  We assume
  $\bm x_i \sim \mbox{Dir}(\alpha_1,\cdots,\alpha_k)$. Define
  $\alpha_0 = \sum_h \alpha_h$ and
  $\bm \alpha = (\alpha_1,\cdots,\alpha_k)^T$.

  We use the standard basis for encoding. We encode $y_{ij} = c_j$
  as $\bm b_{ij} \in \mathbb{R}^{d_j}$ a binary ($0/1$) vector with the $c_j$th
  coordinate being $1$ and all others being $0$. Similarly, we
  encode the membership variable $m_{ij}$ as a $k$ dimensional binary
  vector $\bm m_{ij} \in \mathbb{R}^k$. Consider the first moment of
  $\bm b_{ij}$.
  \begin{align}
    \bm \mu_j = \mathbb{E}(\bm b_{ij}) = \mathbb{E}[\mathbb{E}(\bm
    b_{ij} | \bm m_{ij})] = \mathbb{E}(\bm \Phi_j \bm m_{ij}) =
    \mathbb{E}[\mathbb{E}(\bm \Phi_j \bm m_{ij} | \bm x_i)] =
    \mathbb{E}(\bm \Phi_j \bm x_i) = \bm \Phi_j \frac{\bm \alpha}{\alpha_0}, \notag
  \end{align}
  where $\bm \Phi_j = (\bm \phi_{j1},\cdots \bm \phi_{jk})$. 

  We consider second order moment conditions. There are four types of
  second moments: same variable same subject (type SS), same variable
  cross subject (type SC), cross variable same subject (type CS), and
  cross variable cross subject (type CC). Of the four types,
  only the CS type is needed to prove the theorem.
  The CS type second moment for $\bm b_{ij}$ and $\bm b_{it}$
  ($j \ne t$) can be written as
  \begin{align}
    \mathbb{E}(\bm b_{ij} \circ \bm b_{it}) = \bm \Phi_j \mathbb{E}(\bm m_{ij} \circ \bm m_{it}) \bm
    \Phi_t^T = \bm \Phi_j \mathbb{E}(\bm x_i \circ
    \bm x_i) \bm \Phi_t^T. \notag
  \end{align}
  For a Dirichlet distributed variable,
  \begin{align}
    \mathbb{E}(\bm x_i \circ \bm x_i) &= \cov(\bm x_i) + \mathbb{E}(\bm
    x_i) \circ \mathbb{E}(\bm x_i) \notag \\
    &=\frac{1}{\alpha_0 (\alpha_0 + 1)}\diag(\bm \alpha) +
    \frac{\alpha_0}{\alpha_0^2(\alpha_0+1)} \bm \alpha \circ \bm
    \alpha . \notag
  \end{align}
  Then we have
  \begin{align}
    \mathbb{E}(&\bm b_{ij} \circ \bm b_{it}) = \bm \Phi_j \mathbb{E}(\bm x_i \circ
      \bm x_i) \bm \Phi_t^T \notag\\
    & ~ = \frac{1}{\alpha_0(\alpha_0 + 1)} \sum_{h=1}^k \alpha_h \bm \phi_{jh}
      \circ \bm \phi_{th} + \frac{\alpha_0}{\alpha_0 + 1}\bm \mu_j
      \circ \bm \mu_t. \label{eq:2ndMbjiti}
  \end{align}
  We next consider third order moment conditions.
  There are eight different types of third order moments for
  $\bm b_{ij}$. Only the moments with different variables for the same
  subject are needed to prove the theorem.

  We consider the third cross moment for $\bm b_{ij}$, $\bm b_{it}$
  and $\bm b_{is}$ with $j \ne t \ne s$ for the same subject.
  First we calculate
  $\mathbb{E}(\bm m_{ij} \circ \bm m_{is} \circ \bm m_{it})$.
  \begin{align}
    \mathbb{E}(&\bm m_{ij} \circ \bm m_{is} \circ \bm m_{it}) =
                 \mathbb{E}[\mathbb{E}(\bm m_{ij} \circ \bm m_{is} \circ \bm
                 m_{it} | \bm x_i)] \notag \\
               &=\mathbb{E}(\bm x_i \circ \bm x_i \circ \bm x_i) \notag \\
               &= \frac{1}{\alpha_0(\alpha_0+1)(\alpha_0+2)} \bigg( (\bm \alpha \circ
                 \bm \alpha \circ \bm \alpha) + \sum_{h=1}^k \alpha_h (\bm e_h
                 \circ \bm e_h \circ \bm \alpha ) \notag \\
               & \qquad \qquad +  \sum_{h=1}^k \alpha_h (\bm e_h
                 \circ \bm \alpha \circ \bm e_h ) + \sum_{h=1}^k \alpha_h (\bm
                 \alpha \circ \bm e_h \circ \bm e_h ) \notag \\
               & \qquad \qquad + 2\sum_{h = 1}^k \alpha_h (\bm e_h \circ \bm e_h
                 \circ \bm e_h)\bigg). \notag
  \end{align}
  Here $\bm e_h$ is standard basis vector of length $k$ with $h$th
  coordinate being one.  The third order moment tensor of
  $\bm b_{ij} \circ \bm b_{is} \circ \bm b_{it}$ can be derived as
  \begin{align}
    \mathbb{E}(&\bm b_{ij} \circ \bm b_{is} \circ \bm b_{it})
    =
      \mathbb{E}(\bm m_{ij} \circ \bm m_{is} \circ \bm m_{it} ) \times
      \{\bm \Phi_j, \bm \Phi_s, \bm \Phi_t\} \notag \\
    &= \frac{1}{\alpha_0 (\alpha_0+1) (\alpha_0+2)} \bigg((\bm \Phi_j
      \bm \alpha) \circ (\bm \Phi_s \bm \alpha) \circ (\bm \Phi_t
      \bm \alpha) + \sum_{h=1}^k \alpha_h [ \bm \phi_{jh} \circ \bm
      \phi_{sh} \circ (\bm \Phi_t \bm \alpha )] \notag \\
    & \qquad \qquad + \sum_{h=1}^k \alpha_h [ \bm \phi_{jh} \circ (\bm \Phi_s \bm
      \alpha) \circ \bm \phi_{th}] + \sum_{h=1}^k \alpha_h [ (\bm \Phi_j
      \bm \alpha) \circ \bm \phi_{sh} \circ \bm \phi_{th} ] \notag
    \\
    & \qquad \qquad + 2 \sum_{h=1}^k \alpha_h \bm \phi_{jh} \circ \bm \phi_{sh}
      \circ \bm \phi_{th} \bigg) \notag \\
    &= \frac{1}{\alpha_0 (\alpha_0+1)(\alpha_0+2)} \bigg(\alpha_0^3 \bm
      \mu_j \circ \bm \mu_s \circ \bm \mu_t \notag \\
    & \qquad \qquad + \alpha_0^2(\alpha_0 + 1) \mathbb{E}(\bm b_{ij} \circ
      \bm b_{is} \circ \bm \mu_t ) - \alpha_0^3 \bm \mu_j \circ \bm
      \mu_s \circ \bm \mu_t \notag \\
    & \qquad \qquad + \alpha_0^2(\alpha_0 + 1) \mathbb{E}(\bm \mu_j
      \circ \bm b_{is} \circ \bm b_{it} ) - \alpha_0^3 \bm \mu_j
      \circ \bm \mu_s \circ \bm \mu_t \notag \\
    & \qquad \qquad + \alpha_0^2 (\alpha_0 + 1) \mathbb{E}(\bm b_{ij}
      \circ \bm \mu_s \circ \bm b_{it} ) - \alpha_0^3 \bm \mu_j
      \circ \bm \mu_s \circ \bm \mu_t \notag \\
    & \qquad \qquad + 2\sum_{h=1}^k \alpha_h
      \bm \phi_{jh} \circ \bm \phi_{sh} \circ \bm \phi_{th} \bigg). \label{eq:3rdMbjisitiTensor}
  \end{align}

  The theorem follows Equations \eqref{eq:2ndMbjiti} and
  \eqref{eq:3rdMbjisitiTensor}. For non-categorical
  data, we let $\bm b_{ij} \equiv y_{ij}$ and $\bm \phi_{jh}$ is a
  scalar mean parameter for $y_{ij}$. Equations \ref{eq:2ndMbjiti} and
  \ref{eq:3rdMbjisitiTensor} still hold.
\end{proof}

\subsection{Proof of Theorem \ref{consistencyTheorem}}
\label{app:theorem2}
\begin{proof}
  For notation simplicity we suppress $\sups{\bm A}{\cdot}_n$ in
  $\sups{Q}{\cdot}_n(\bm \Phi;\sups{\bm A}{\cdot}_n)$ and
  $\sups{\bm A}{\cdot}$ in
  $\sups{Q}{\cdot}_0(\bm \Phi;\sups{\bm A}{\cdot})$. Lemma
  \ref{uniformLemma} implies
  $\lim_{n\rightarrow \infty}
  \pr[|\sups{Q}{\cdot}_n(\sups{\widehat{\bm \Phi}}{\cdot}) -
  \sups{Q}{\cdot}_0(\sups{\widehat{\bm \Phi}}{\cdot})| < \epsilon/3] =
  1$
  and
  $\lim_{n\rightarrow \infty} \pr[|\sups{Q}{\cdot}_n(\bm \Phi_0) -
  \sups{Q}{\cdot}_0(\bm \Phi_0)| < \epsilon/3] = 1$
  for $\epsilon > 0$. This result also implies
  \begin{align}
    \lim_{n\rightarrow \infty} &\pr[\sups{Q}{\cdot}_0(\sups{\widehat{\bm
    \Phi}}{\cdot}) < \sups{Q}{\cdot}_n(\sups{\widehat{\bm
    \Phi}}{\cdot})+ \epsilon/3] = 1. \label{eq:consistproof1} \\
    \lim_{n\rightarrow \infty} &\pr[\sups{Q}{\cdot}_n(\bm
    \Phi_0) < \sups{Q}{\cdot}_0(\bm
    \Phi_0)+ \epsilon/3] = 1. \label{eq:consistproof2}
  \end{align}
  On the other hand, $\sups{\widehat{\bm \Phi}}{\cdot}$ minimizes
  $\sups{Q}{\cdot}_n(\bm \Phi)$, therefore
  \begin{align}
    \lim_{n\rightarrow \infty} \pr[\sups{Q}{\cdot}_n(\sups{\widehat{\bm
    \Phi}}{\cdot}) < \sups{Q}{\cdot}_n(\bm
    \Phi_0)+ \epsilon/3] = 1. \label{eq:consistproof3}
  \end{align}
  Equations \eqref{eq:consistproof1} and \eqref{eq:consistproof3} imply
  \begin{align}
    \lim_{n\rightarrow \infty} \pr[\sups{Q}{\cdot}_0(\sups{\widehat{\bm
    \Phi}}{\cdot}) < \sups{Q}{\cdot}_n(\bm
    \Phi_0)+ 2\epsilon/3] = 1. \notag
  \end{align}
  Together with Equation \eqref{eq:consistproof2}, we get
  \begin{align}
    \lim_{n\rightarrow \infty} \pr[\sups{Q}{\cdot}_0(\sups{\widehat{\bm
    \Phi}}{\cdot}) < \sups{Q}{\cdot}_0(\bm
    \Phi_0)+ \epsilon] = 1. \notag
  \end{align}
  Therefore
  \begin{align}
    \lim_{n\rightarrow \infty} \pr[0 \le \sups{Q}{\cdot}_0(\sups{\widehat{\bm
    \Phi}}{\cdot}) <  \epsilon] = 1 \label{eq:consistproof4}
  \end{align}
  follows with $\sups{Q}{\cdot}_0(\bm \Phi_0) = 0$ and
  $\sups{Q}{\cdot}_0(\sups{\widehat{\bm \Phi}}{\cdot}) \ge 0$.  Next,
  we choose a neighborhood ${N}$, which contains $\bm \Phi_0$
  in $\bm \Theta$. Due to the compactness of $\bm \Theta$ the
  neighborhood ${ N}^C$ is also compact. The continuousness of
  $\sups{Q}{\cdot}_0(\bm \Phi)$ implies the existence of
  $\inf_{\bm \Phi \in N^C} \sups{Q}{\cdot}_0(\bm \Phi)$ and it is
  positive. Let
  $\epsilon = \inf_{\bm \Phi \in {N}^C} \sups{Q}{\cdot}_0(\bm
  \Phi)$, then we get
  \begin{align}
    \lim_{n\rightarrow \infty} \pr[0 \le \sups{Q}{\cdot}_0(\sups{\widehat{\bm
    \Phi}}{\cdot}) <  \inf_{\bm \Phi \in N^C}
    \sups{Q}{\cdot}_0(\bm \Phi)] = 1. \label{eq:consistproof5}
  \end{align}
  Therefore
  $\lim_{n\rightarrow \infty} \pr(\sups{\widehat{\bm \Phi}}{\cdot}
  \not\in {N}^C) = 1$,
  which suggests
  $\lim_{n\rightarrow \infty} \pr(\sups{\widehat{\bm \Phi}}{\cdot} \in
  {N}) = 1$. Shrinking the neighborhood size of ${N}$ we get
  \begin{align}
    \lim_{n\rightarrow \infty} \pr(\sups{\widehat{\bm \Phi}}{\cdot} =
    \bm \Phi_0 ) = 1. \notag
  \end{align}
\end{proof}

\subsection{Proof of Theorem \ref{normalityTheorem}}
\label{app:theorem3}
\begin{proof}
  We approximate
  $\sups{\bm f}{\cdot}_n(\sups{\widehat{\bm \Phi}}{\cdot})$
  using first order Taylor expansion
  \begin{align}
    \sups{\bm f}{\cdot}_n(\sups{\widehat{\bm \Phi}}{\cdot})
    = \sups{\bm f}{\cdot}_n(\bm \Phi_0) + \sups{\bm
    G}{\cdot}_n(\bm \Phi_0)[\mbox{vec}(\sups{\widehat{\bm
    \Phi}}{\cdot}) - \mbox{vec}(\bm \Phi_0)] + O\{[\mbox{vec} (
    \sups{\widehat{\bm \Phi}}{\cdot}) - \mbox{vec}(\bm \Phi_0)]^2\}
    \notag
  \end{align}
  Ignoring the high order term, we left multiply both sides by
  $[\sups{\bm G}{\cdot}_n(\sups{\widehat{\bm \Phi}}{\cdot})]^T
  \sups{\bm A}{\cdot}_n$. Then we get
  \begin{align}
    [\sups{\bm G}{\cdot}_n(\sups{\widehat{\bm \Phi}}{\cdot})]^T
  \sups{\bm A}{\cdot}_n &\sups{\bm
    f}{\cdot}_n(\sups{\widehat{\bm \Phi}}{\cdot}) \approx [\sups{\bm
    G}{\cdot}_n(\sups{\widehat{\bm \Phi}}{\cdot})]^T \sups{\bm
    A}{\cdot}_n  \sups{\bm f}{\cdot}_n(\bm \Phi_0) \notag \\
    &+ [\sups{\bm G}{\cdot}_n(\sups{\widehat{\bm \Phi}}{\cdot})]^T
    \sups{\bm A}{\cdot}_n \sups{\bm  G}{\cdot}_n(\bm
    \Phi_0)[\mbox{vec}(\sups{\widehat{\bm  \Phi}}{\cdot}) -
    \mbox{vec}(\bm \Phi_0)]. \notag
  \end{align}
  The fact that estimator $\sups{\widehat{\bm \Phi}}{\cdot}$ minimizes
  $\sups{Q}{\cdot}_n(\bm \Phi, \sups{\bm A}{\cdot}_n)$ implies the
  left hand side equals to zero. Therefore we get
  \begin{align}
    n^{1/2}[\mbox{vec}(\sups{\widehat{\bm
    \Phi}}{\cdot}) - \mbox{vec}(\bm \Phi_0)] \approx - \{[\sups{\bm
    G}{\cdot}_n(\sups{\widehat{\bm \Phi}}{\cdot})]^T \sups{\bm
    A}{\cdot}_n \sups{\bm
    G}{\cdot}_n(\bm \Phi_0) \}^{-1} [\sups{\bm
    G}{\cdot}_n(\sups{\widehat{\bm \Phi}}{\cdot})]^T \sups{\bm
    A}{\cdot}_n n^{1/2} \sups{\bm f}{\cdot}_n(\bm
    \Phi_0). \notag
  \end{align}
  The theorem follows with $n^{1/2} \sups{\bm f}{\cdot}_n(\bm
    \Phi_0) \overset{p}{\rightarrow} \N(\bm 0,\sups{\bm S}{\cdot})$ and
    Assumptions \ref{assumption1} and \ref{assumption2}.
\end{proof}

\subsection{Derivatives of moment functions}
\label{app:derivatives}
\subsubsection{Second moment matrix}
The second moment matrix $\sups{\bm F}{2}_{jt}(\bm y_i, \bm \Phi)$ may be written as $\bm b_{ij} \circ \bm b_{it} - \bm
\Phi_j \E(\bm x_i \circ \bm x_i) \bm \Phi_t^T$.  The derivatives of
$\sups{\bm F}{2}_{jt}(\bm y_i, \bm \Phi)$ with respect to $\bm \Phi_j$
and $\bm \Phi_t$ can be written as
\begin{align}
  \frac{\partial \mbox{vec}[\sups{\bm F}{2}_{jt}(\bm y_i, \bm
  \Phi)]}{\partial \mbox{vec}(\bm \Phi_j)} &= - \frac{\partial\{ [\bm
                                             \Phi_t \E(\bm x_i 
  \circ \bm x_i)] \otimes \bm I_{d_j}\}\mbox{vec}(\bm
  \Phi_j)}{\partial \mbox{vec}(\bm \Phi_j)} \notag \\
  &= - [\bm \Phi_t \E(\bm x_i
  \circ \bm x_i)] \otimes \bm I_{d_j}, \notag \\
  \frac{\partial \mbox{vec}[\sups{\bm F}{2}_{jt}(\bm y_i, \bm
  \Phi)]}{\partial \mbox{vec}(\bm \Phi_t)}
  &= - \bm T \frac{\partial\mbox{vec}[\bm \Phi_t \E(\bm x_i \circ
    \bm x_i) \bm \Phi_j^T]}{\partial \mbox{vec}(\bm \Phi_t)} \notag \\
  &= - \bm T \frac{\partial \{[\bm \Phi_j \E(\bm x_i \circ \bm x_i)] \otimes \bm
    I_{d_t} \mbox{vec}(\bm \Phi_t)\} }{\partial \mbox{vec}(\bm
    \Phi_t)} \notag \\
  &= - \bm T \{ [\bm \Phi_j \E(\bm x_i \circ \bm x_i)] \otimes \bm
    I_{d_t} \}, \notag
\end{align}
where $\otimes$ indicates a Kronecker product and $\bm T$ is a
$d_tk \times d_tk$ 0/1 matrix that satisfies
\begin{align}
  \mbox{vec}[\bm \Phi_j \E(\bm x_i \circ \bm x_i) \bm \Phi_t^T] =
  \bm T \mbox{vec}[\bm \Phi_t \E(\bm x_i \circ \bm x_i) \bm
  \Phi_j^T]. \notag
\end{align}
Therefore
$\E[\partial \sups{\bm f}{2}(\bm y_i, \bm \Phi)/\partial \bm \Phi]$ is
a block matrix with block of
$- \bm \Phi_t^T \E(\bm x_i \circ \bm x_i) \otimes \bm I_{d_j}$ on
columns corresponding to $\mbox{vec}(\bm \Phi_t)$ and rows
corresponding to
$\mbox{vec}[\sups{\bm F}{2}_{jt}(\bm y_i, \bm \Phi)]$.

\subsubsection{Third moment tensor}

We next consider the third moment tensor. Write
$\sups{\bm F}{3}_{jst}(\bm y_i,\bm \Phi)$ as
$\bm b_{ij} \circ \bm b_{ij} \circ \bm b_{ij} - \E(\bm x_i \circ \bm
x_i \circ \bm x_i) \times_1 \bm \Phi_j \times_2 \bm \Phi_s \times_3
\bm \Phi_t$. Then only the second term involves $\bm \Phi$.

The derivatives of $\sups{\bm F}{3}_{jst}(\bm y_i,\bm \Phi)$ with
respect to $\bm \Phi_j$ can be written as
\begin{align}
  \frac{\partial \mbox{vec}[\sups{\bm F}{3}_{jst}(\bm y_i,\bm
  \Phi)]}{\partial \mbox{vec}(\bm \Phi_j)}
  &=  - \frac{\partial \mbox{vec}\{[\E(\bm x_i \circ \bm x_i \circ \bm
    x_i)  \times_1 \bm \Phi_j \times_2 \bm \Phi_s \times_3 \bm
    \Phi_t]_{(1)}\}} {\partial \mbox{vec}(\bm \Phi_j)} \notag \\
  & = -\frac{\partial \mbox{vec}\{[\bm \Phi_j \E(\bm x_i \circ \bm x_i
    \circ \bm x_i)_{(1)} ( \bm \Phi_t \otimes \bm \Phi_s )^T\}}{ \partial \mbox{vec}(\bm \Phi_j)} \notag \\
   & = - (\bm \Phi_t \otimes \bm \Phi_s) \mbox{vec}[ \E(\bm x_i  \circ
     \bm x_i \circ  \bm x_i)_{(1)}]^T \otimes \bm I_{d_j}, \notag
\end{align}
where subscript $(1)$ indicates model-1 unfolding of a three way
tensor.

The derivatives of $\sups{\bm F}{3}_{jst}(\bm y_i,\bm \Phi)$ with
respect to $\bm \Phi_s$ and $\bm \Phi_t$ can be calculated accordingly
by introducing 0/1 transformation matrices $\bm T_{(2)(1)}$ and
$\bm T_{(3)(1)}$ both with size $d_jd_sd_t \times d_jd_sd_t$ that
satisfy
\begin{align}
  \mbox{vec}\{[\sups{\bm F}{3}_{jst}(\bm y_i,\bm \Phi)]_{(1)}\}
  &=  \bm T_{(2)(1)} \mbox{vec}\{ [\sups{\bm F}{3}_{jst}(\bm y_i,\bm \Phi)]_{(2)}\} \notag \\
  &= \bm T_{(3)(1)} \mbox{vec}\{ [\sups{\bm F}{3}_{jst}(\bm y_i,\bm \Phi)]_{(3)}\}. \notag
\end{align}

Then
\begin{align}
  \frac{\partial \mbox{vec}[\sups{\bm F}{3}_{jst}(\bm y_i,\bm
  \Phi)]}{\partial \mbox{vec}(\bm \Phi_s)}
  &= - \bm T_{(2)(1)} \frac{\partial \mbox{vec}\{[\E(\bm x_i \circ \bm x_i \circ \bm
    x_i)  \times_1 \bm \Phi_j \times_2 \bm \Phi_s \times_3 \bm
    \Phi_t]_{(2)}\}} {\partial \mbox{vec}(\bm \Phi_s)} \notag \\
  & = -  \bm T_{(2)(1)} \frac{\partial \mbox{vec}\{[\bm \Phi_s \E(\bm x_i \circ \bm x_i
    \circ \bm x_i)_{(2)} ( \bm \Phi_t \otimes \bm \Phi_j )^T\}}
    {\partial \mbox{vec}(\bm \Phi_s)}  \notag \\
  & = - \bm T_{(2)(1)} \{ (\bm \Phi_t \otimes \bm \Phi_j) [ \E(\bm x_i \circ
    \bm x_i \circ \bm x_i)_{(2)}]^T \otimes \bm I_{d_s}\} , \notag\\
  \frac{\partial \mbox{vec}[\sups{\bm F}{3}_{jst}(\bm y_i,\bm \Phi)]}
  {\partial \mbox{vec}(\bm \Phi_t)}
  &= -  \bm T_{(3)(1)} \frac{\partial \mbox{vec}\{[\E(\bm x_i \circ \bm x_i \circ \bm
    x_i)  \times_1 \bm \Phi_j \times_2 \bm \Phi_s \times_3 \bm
    \Phi_t]_{(3)}\}} {\partial \mbox{vec}(\bm \Phi_t)} \notag \\
  & = - \bm T_{(3)(1)} \frac{\partial \mbox{vec}\{[\bm \Phi_t \E(\bm x_i \circ \bm x_i
    \circ \bm x_i)_{(3)} ( \bm \Phi_s \otimes \bm \Phi_j
    )^T\}}{\partial \mbox{vec}(\bm \Phi_t)}  \notag \\
  & = -  \bm T_{(3)(1)} \{ (\bm \Phi_s \otimes \bm \Phi_j) [ \E(\bm x_i \circ
    \bm x_i \circ \bm x_i)_{(3)}]^T \otimes \bm I_{d_t} \} . \notag
\end{align}

The conditions 1), 3) and 4) in Assumption 2 in main text follow after
we calculating the derivatives of moment functions.

\subsection{Derivation of Newton-Raphson update}
\label{app:NRupdate}

We denote
\begin{align}
  \sups{\bm E}{2}_{n,jt}
  &= \sups{\bm F}{2}_{n,jt}(\bm \Phi) + \bm \Phi_j
    \sups{\bm \Lambda}{2} \bm \Phi_t^T, \notag \\
  \sups{\bm E}{3}_{n,jst}
  & = \sups{\bm F}{3}_{n,jst}(\bm \Phi) + \sups{\bm
    \Lambda}{3} \times_1 \bm \Phi_j \times_2 \bm
    \Phi_s \times_3 \bm
    \Phi_t. \notag
\end{align}
With identity matrix, then the two objective functions can be written
as
\begin{align}
  \sups{Q}{2}_n(\bm \Phi) &= \sum_{j=1}^{p-1} \sum_{t = j+1}^p
                                 ||\sups{\bm
  E}{2}_{n,jt}  - \bm \Phi_j \sups{\bm \Lambda}{2} \bm
  \Phi_t^T||^2_F,  \notag \\
  \sups{Q}{3}_n(\bm \Phi) &= \sum_{j=1}^{p-1} \sum_{t = j+1}^p
                                 ||\sups{\bm
  E}{2}_{n,jt}  - \bm \Phi_j \sups{\bm \Lambda}{2} \bm
  \Phi_t^T||^2_F \notag \\
    & \qquad +
  \sum_{j=1}^{p-2} \sum_{s = j+1}^{p-1} \sum^{p}_{t = s+1}
  ||\sups{\bm E}{3}_{n,jst} - \sups{\bm
  \Lambda}{3} \times_1 \bm \Phi_j \times_2 \bm
  \Phi_s \times_3 \bm \Phi_t \notag||_F^2. \notag
\end{align}
Here we suppress the weight matrix in the objective functions.  We
first consider $\sups{Q}{2}_n(\bm \Phi)$. The terms involve
$\bm \phi_{jh}$ are
\begin{align}
  \sum_{t=1,t\ne j}^p \bigg[-2\sups{\lambda}{2}_h (\sups{\overline{\bm
  E}}{2}_{n,jt} \bm \phi_{th})^T \bm \phi_{jh} +
  (\sups{\lambda}{2}_h)^2 (\bm
  \phi_{th}^T \bm \phi_{th} ) \bm \phi_{jh}^T \bm \phi_{jh}\bigg], \notag
\end{align}
where
$\sups{\overline{\bm E}}{2}_{n,jt} = \sups{\bm E}{2}_{n,jt} - \sum_{h'
  \ne h} \sups{\lambda}{2}_{h'} \bm \phi_{jh'} \circ \bm \phi_{th'}$
and $\sups{\lambda}{2}_h$ is the $h$th diagonal element of
$\sups{\bm \Lambda}{2}$. Here we use the fact that
$||\sups{\bm E}{2}_{n,jt} - \bm \Phi_j \sups{\bm \Lambda}{2} \bm
\Phi_t^T||^2_F = ||\sups{\bm E}{2}_{n,tj} - \bm \Phi_t \sups{\bm
  \Lambda}{2} \bm \Phi_j^T||^2_F$. By letting
\begin{align}
  \sups{\bm \xi}{2} &= -2 \sups{\lambda}{2}_h \sum_{t=1, t\ne j}^p
                      (\sups{\overline{\bm
  E}}{2}_{n,jt} \bm \phi_{th} ), \notag \\
  \sups{\gamma}{2} &= (\sups{\lambda}{2}_h)^{2} \sum_{t=1,t\ne j}^p \bm
  \phi_{th}^T \bm \phi_{th}, \notag
\end{align}
the gradient $\nabla \sups{Q}{2}_n(\bm \phi_{jh})$ and Hessian
$\nabla^2 \sups{Q}{2}_n(\bm \phi_{jh})$ can be written as
\begin{align}
  \nabla \sups{Q}{2}_n(\bm \phi_{jh}) &= \sups{\bm \xi}{2} +
  2\sups{\gamma}{2} \bm \phi_{jh}, \notag\\
  \nabla^2 \sups{Q}{2}_n(\bm \phi_{jh}) &= 2 \sups{\gamma}{2} \bm I. \notag
\end{align}
The update rule in \eqref{eq:updatePhiM2} can be derived accordingly.

Then we consider $\sups{Q}{3}_n(\bm \Phi)$. The terms involve $\bm
\phi_{jh}$ are
\begin{align}
  & \sum_{t=1,t \ne j}^p \bigg[-2
    \sups{\lambda}{2}_h  (\sups{\overline{\bm
    E}}{2}_{n,jt} \bm \phi_{th})^T \bm \phi_{jh} +
    (\sups{\lambda}{2}_h)^2 (\bm
    \phi_{th}^T \bm \phi_{th} ) \bm \phi_{jh}^T \bm \phi_{jh}
    \bigg] \notag \\
  & \qquad + \sum_{s=1,s\ne j}^p \sum_{t=1,t \ne s,t
    \ne j}^p \bigg[ -2  \la \sups{\overline{\bm
    E}}{3}_{n,jst},  \sups{\lambda}{3}_h \bm
    \phi_{jh} \circ \bm \phi_{sh} \circ \bm \phi_{th} \ra +
    ||\sups{\lambda}{3}_h \bm \phi_{jh} \circ \bm \phi_{sh} \circ
    \bm \phi_{th}||^2_F \bigg], \notag
\end{align}
where
$\sups{\overline{\mb E}}{2}_{n,jt} = \sups{\bm E}{2}_{n,jt} -
\sum_{h' \ne h} \sups{\lambda}{2}_{h'} \bm \phi_{jh'} \circ \bm
\phi_{th'}$
and
$\sups{\overline{\bm E}}{3}_{n,jst} = \sups{\bm E}{3}_{n,jst} -
\sum_{h' \ne h} \sups{\lambda}{3}_{h'} \bm \phi_{jh'} \circ \bm
\phi_{sh'} \circ \bm \phi_{th'}$.
Again we use the symmetric property of
$||\sups{\bm E}{2}_{n,jt} - \bm \Phi_j \sups{\bm \Lambda}{2}
\bm \Phi_t^T ||^2_F$
and
$||\sups{\bm E}{3}_{n,jst} - \sups{\bm \Lambda}{3} \times_1
\bm \Phi_j \times_2 \bm \Phi_s \times_3 \bm \Phi_t ||_F^2$.
By organizing the terms, we get
\begin{align}
  \sum_{t=1,t \ne j}^p &\bigg[-2
    \sups{\lambda}{2}_h  (\sups{\overline{\bm
  E}}{2}_{n,jt} \bm \phi_{th})^T \bm \phi_{jh} + (\sups{\lambda}{2}_h)^2  (\bm
  \phi_{th}^T \bm \phi_{th} ) \bm \phi_{jh}^T \bm \phi_{jh}\bigg]
  \notag \\
  & \qquad  + \sum_{s=1,s \ne j}^p\sum_{t=1, t \ne s, t
    \ne j}^p\bigg[ -2 \sups{\lambda}{3}_h (\sups{\overline{\bm E}}{3}_{n,jst}
    \times_2 \bm \phi_{sh} \times_3 \bm \phi_{th})^T \bm \phi_{jh} +
    (\sups{\lambda}{3}_h)^2 (\bm \phi^T_{sh} \bm \phi_{sh}) (\bm
    \phi^T_{th} \bm \phi_{th}) (\bm \phi^T_{jh} \bm \phi_{jh}) \bigg]. \notag
\end{align}
By letting
\begin{align}
  \sups{\bm \xi}{3} &= -2\sups{\lambda}{2}_h  \sum_{t=1,t\ne j}^p
                      (\sups{\overline{\bm E}}{2}_{n,jt} \bm \phi_{th})
 -2 \sups{\lambda}{3}_h \sum_{s=1,s\ne j}^p \bigg[
    \sum_{t=1,t\ne s, t\ne j}^p (\sups{\overline{\bm E}}{3}_{n,jst}
    \times_2 \bm \phi_{sh} \times_3 \bm \phi_{th})\bigg], \notag \\
  \sups{\gamma}{3} &= (\sups{\lambda}{2}_h)^2 \sum_{t=1,t\ne
                j}^p\bm \phi_{th}^T \bm \phi_{th}
 +  (\sups{\lambda}{3}_h)^2  \sum_{s=1,s\ne j}^p
    \bigg[ \sum_{ t=1,t\ne s, t \ne j}^p (\bm \phi^T_{sh} \bm
                     \phi_{sh}) (\bm \phi^T_{th} \bm \phi_{th})\bigg],
                     \notag
\end{align}
the gradient $\nabla \sups{Q}{3}_n(\bm \phi_{jh})$ and Hessian
$\nabla^2 \sups{Q}{3}_n(\bm \phi_{jh})$ can be written as
\begin{align}
  \nabla \sups{Q}{3}_n(\bm \phi_{jh}) &= \sups{\bm \xi}{3} +
  2\sups{\gamma}{3} \bm \phi_{jh}, \notag\\
  \nabla^2 \sups{Q}{3}_n(\bm \phi_{jh}) &= 2 \sups{\gamma}{3} \bm I. \notag
\end{align}
The update rule in \eqref{eq:updatePhiM3} follows directly.

\bigskip


\section*{References}
\bibliographystyle{apalike}
\bibliography{References}
\end{document}